\documentclass[11pt]{amsart}

\title{Cramér distance and discretizations of circle expanding maps I: theory}
\date{\today}

\author{Pierre-Antoine Guih\'eneuf}
\address{Pierre-Antoine Guih\'eneuf: Sorbonne Université and Université de Paris, CNRS, IMJ-PRG, F-75005 Paris, France.}
\email{pierre-antoine.guiheneuf@imj-prg.fr}

\author{Maurizio Monge}
\address{Maurizio Monge: Instituto de Matem{\`a}tica da UFRJ, Av. Athos da Silveira Ramos 149, Centro de Tecnologia, Bloco C Cidade Univesit{\`a}ria, Ilha do Fund{\~a}o, Caixa Postal 68530 21941-909, Rio de Janeiro, RJ, Brasil}
\email{maurizio.monge@im.ufrj.br}


\usepackage[utf8]{inputenc}
\usepackage[english]{babel}
\usepackage[T1]{fontenc}
\usepackage{amsfonts}
\usepackage{dsfont}
\usepackage{amsmath}
\usepackage{amssymb}
\usepackage{stmaryrd}
\usepackage{amsthm}
\usepackage{enumerate}
\usepackage{pgf,tikz}
\usetikzlibrary{decorations.pathreplacing, shapes.multipart, arrows, matrix}
\usepackage[pdftex,colorlinks=true,linkcolor=blue,citecolor=blue,urlcolor=blue]{hyperref}

\usepackage{mathtools}

\newtheorem{lemme}{Lemma}[section]

\newtheorem{theoreme}[lemme]{Theorem}
\newtheorem{prop}[lemme]{Proposition}
\newtheorem{coro}[lemme]{Corollary}

\newtheorem{add}[lemme]{Addendum}

\newtheorem*{ques}{Question}
\newtheorem{theo}{Theorem}

\newtheorem*{theo*}{Theorem}

\theoremstyle{definition}
\newtheorem{definition}[lemme]{Definition}

\theoremstyle{remark}
\newtheorem{rem}[lemme]{Remark}

\numberwithin{equation}{section}

\newcommand{\Er}{\mathcal{E}}

\newcommand{\N}{\mathds{N}}

\newcommand{\R}{\mathds{R}}

\newcommand{\T}{\mathds{T}}

\newcommand{\Q}{\mathds{Q}}
\newcommand{\Z}{\mathds{Z}}
\newcommand{\D}{\mathcal{D}}

\newcommand{\varep}{\varepsilon}

\newcommand{\Leb}{\operatorname{Leb}}

\newcommand{\ud}{\,\mathrm{d}}
\newcommand{\Prb}{\mathcal{P}}

\newcommand{\card}{\operatorname{Card}}
\newcommand{\dist}{\operatorname{dist}}

\newcommand{\1}{\mathds 1}

\newcommand{\disc}{\operatorname{d_C}}

\newcommand{\Sp}{\mathds{S}}

\newcommand{\ind}{{\boldsymbol{i}}}
\newcommand{\Disc}{\operatorname{d_C}}
\newcommand{\len}{\operatorname{length}}
\newcommand{\fat}{\operatorname{\wp}}
\newcommand{\Ll}{L_f}

\newcommand{\pr}{\operatorname{pr}}
\newcommand{\Corr}{\operatorname{Corr}}
\newcommand{\Ell}{\ensuremath{\boldsymbol\ell}}
\newcommand{\Varep}{\ensuremath{\boldsymbol\varepsilon}}
\newcommand{\Var}{\operatorname{Var}}

\allowdisplaybreaks

\begin{document}

\begin{abstract}
This paper is aimed to study the ergodic short-term behaviour of discretizations of circle expanding maps. More precisely, we prove some asymptotics of the distance between the $t$-th iterate of Lebesgue measure by the dynamics $f$ and the $t$-th iterate of the uniform measure on the grid of order $N$ by the discretization on this grid, when $t$ is fixed and the order $N$ goes to infinity. This is done under some explicit genericity hypotheses on the dynamics, and the distance between measures is measured by the mean of \emph{Cramér} distance. The proof is based on a study of the corresponding linearized problem, where the problem is translated into terms of equirepartition on tori of dimension exponential in $t$.

A numerical study associated to this work is presented in \cite{paper2}.
\end{abstract}

\maketitle

\setcounter{tocdepth}{1}
\tableofcontents

\section{Introduction}

\subsection*{Motivations}

In one of the last papers he published \cite{MR1678095}, Oscar E. Lanford was proposing to study the behaviour of spatial discretizations of expanding maps of the circle in some limiting regime. The question was to decide whether in most of cases, the middle-term ergodic behaviour of such discretizations reflects the actual dynamics of the map.

To fix the notations, let us take $f : \Sp^1\to \Sp^1$ an expanding map (meaning that $f'(x)>1$ for any $x \in\Sp^1$) and consider the grid $E_N$ made of $N$ points of $\Sp^1$ equally spaced. That is, identifying $\Sp^1\simeq \R/\Z$ with $[0,1]$, one sets
\[E_N = \left\{\frac{i}{N} \ \big|\  0\le i \le N-1 \right\}.\]
To each of these grids is associated a projection $P_N : \Sp^1\to E_N$ on the nearest point of $E_N$ (for some points one has two choices for the nearest neighbour, one does a choice once for all and this choice will not play any role in the sequel). This leads to the definition of the \emph{discretization} of the map $f$ relatively to the grid $E_N$ as
\begin{align*}
f_N : E_N & \longrightarrow E_N\\
x & \longmapsto P_N(f(x)).
\end{align*}
Remark that if $N = 2^n$, this corresponds to a discretization realizing the rounding of $f(x)$ with $n$ binary digits.

In \cite{MR1678095}, Lanford asked whether the dynamics of the maps $f_N^k$ looks like the one of $f^k$ in the regime $\log N \ll k \ll \sqrt N$. Here is how he justifies these bounds: \textit{``The first $\ll$ allows the computed orbit to deviate macroscopically from the true one over most of its length; the second is in any case usually satisfied in practice and ought to mean that the times considered are short enough so that the effects of the strict finiteness of the space of states are not important. In fact: it might be prudent to replace the second $\ll$ by the stronger condition $\log k \ll \log N$.''} His article includes enlightening philosophical thoughts, supported by some numerical experiments.

This research program was tackled by Paul P.~Flockermann in his unpublished PhD thesis \cite{Flocker} (under the supervision of Lanford). In this work he obtains partial results towards the limiting behaviour of the ``non-injectivity'' of the maps $f_N$ -- i.e. the quantity $\card(f_N^k(E_N))/\card(E_N)$ -- when $k$ is fixed and $N$ goes to infinity\footnote{More precisely, he treats the corresponding linear case: he gets some local statements.}. These statements are valid under genericity assumptions on the expanding map\footnote{They become false for the trivial example of $x\mapsto 2x$.} $f$: they concern either generic $C^r$ expanding maps, for $1<r\le +\infty$, or any real-analytic expanding map different from $x\mapsto 2x$. Lanford and Flockermann were writing an article to complete these partial results (they had an unpublished draft) which unfortunately has never been published.

These results had been obtained independently by Vladimirov in \cite{vladimirov2015quantized} (further works based on this grounding article were published in \cite{MR1661120, MR1769577, MR1956409, MR1894464}). In this article, the author founds a solid theoretical basis about the discretizations' behaviour:
the algebras of quasiperiodic subsets of the lattice, their statistical properties (frequency measurability) under nonresonance conditions, their algebraic properties with respect to discretizations of linear maps, the role of skew products of measure-preserving automorphisms of multidimensional tori in the asymptotic independence and uniform distribution of quantisation errors, etc.
This approach reveals more powerful than Flockermann's approach: in addition to the equidistribution of roundoff errors, Vladimirov gets the fact that the asymptotic rate of injectivity is 0, and some functional central limit theorem, which was published with Vivaldi in \cite{MR2031150}. Early apparitions of this kind of ideas can be found in the work of Voevodin \cite{MR219227}.

One of these results has been re-discovered independently by the first author in \cite{Gui15a}: it is proved that the actual limit of the non injectivity rate $\card(f_N^k(E_N))/\card(E_N)$ when $k$ is fixed and $N$ goes to infinity is zero. The techniques used in this article are a bit different from the ones of Flockermann and Vladimirov: they involve the notion of ``model set'', usually used in the study of quasicrystals, and some theorems from basic geometry of numbers, and allow to get similar results in different settings. The approach of the present article is based on these techniques.

\subsection*{Main results}

The aim of this article is to make a contribution in the direction of Lanford's program, by looking at the short-term ergodic behaviour of discretizations: we will compare the actions of the maps $f$ and $f_N$ on respectively Lebesgue measure $\Leb$ on $\Sp^1$ and the uniform measure $\Leb_N$ on $E_N$. This comparison will be made using a distance on measures called here \emph{Cramér distance} and denoted $\Disc$, which spans the weak-* topology (see Section~\ref{SecDistMeas}). This distance il also called Cramér-von Mises distance, or ``the $L_2$-metrics between distribution functions'' \cite{MR1105086, zbMATH05216873}. Our goal will be to get an asymptotics for the quantity
\begin{equation}\label{EqStudied}
\Disc\big(f^k_*(\Leb), (f_N^k)_*(\Leb_N)\big)
\end{equation}
when $k$ is fixed and $N$ goes to infinity.

As for the previous works already described, we will need genericity assumptions to ensure that there is no phenomenon of resonance between the dynamics and the grid. For example, if $f(x) = 2x \mod 1$ and $N=2^n$, then the orbit of any point of $E_N$ under $f_N$ eventually falls in the fixed point $0$; this is a very specific phenomenon that one wants to avoid to understand what happens ``in most cases''\footnote{Following Lanford in \cite{MR1678095}, it is interesting to note that this hypothesis of genericity, which ensures some uniform repartition properties at a mesoscopic scale (e.g. Proposition \ref{RoundoffLin}), is also a classical assumption for the problem of deriving fluid mechanics laws from a microscopical deterministic model. See for example \cite{MR3455156}.}. Hence, we will consider \emph{generic properties} on the spaces of $C^r$ expanding maps: a property will be called generic if satisfied on at least a countable intersection of open and dense sets of $C^r$ expanding maps. As these sets of maps are Baire, genericity has some nice natural properties: a generic property is satisfied on a dense set of maps, the fact of satisfying two generic properties is generic, etc. In fact, the genericity properties needed to get our results are very weak, so our theorems are also valid under some different genericity assumptions (see \cite{Gui15a} for a discussion).

Our main theorem is the following.

\begin{theo}\label{MainTheo}
Let $r\ge 1$, $f$ a generic $C^r$ expanding map of the circle $\Sp^1$, and $k\in\N$. Then
\begin{equation}\label{EqIntMieuxIntro}
\lim_{N\to +\infty} N^2 \Disc\Big(f^k_*(\Leb), (f_N^k)_*(\Leb_N)\Big)^2
= \frac{1}{12} + \frac{1}{12} \sum_{m=0}^{k-1} \big\langle D(f^{k-m}), (L_f^m 1)^2 \big\rangle,
\end{equation}
where $\langle\cdot,\cdot\rangle$ stands for the $L^2$ scalar product, $L_f$ is the Ruelle-Perron-Frobenius transfer operator defined by \eqref{EqDefRPF}, and $f^{k-m}$ is the $(k-m)$-th iterate of $f$.
\end{theo}

This asymptotics tells us at which speed the measures $f^k_*(\Leb)$ and $(f_N^k)_*(\Leb_N)$ move apart one from the other. More precisely, from this theorem and an estimation of the terms in Lasota-Yorke inequality (see \cite{MR3105868} or \cite[Theorem 12.17]{Guih-These}), one can easily deduce  the following.

\begin{coro}
Let $r > 1$, and $f$ a generic $C^r$ expanding map of $\Sp^1$. Then there exist two constants $1<c<C$, depending only on $\inf f'$ and the $C^{r-1}$ norm of $f'$, such that for any $K\in\N$, there is $N_0\in\N$ such that for any $N\ge N_0$ and any $k\le K$, one has
\[
N^2 \Disc\Big(f^k_*(\Leb), (f_N^k)_*(\Leb_N)\Big)^2 \in [c^k,C^k].
\]
\end{coro}

In other words, for $K$ fixed, if $N$ is large enough, then the behaviour of the Cramér distance will be typically exponential for times smaller than $K$.
\bigskip

Let us describe a bit the term $f^k_*(\Leb)$. As realized quite a long time ago by physicists, and then in the 70's by mathematicians including David Ruelle, the action of a hyperbolic map on measures can be discribed with the help of the \emph{Ruelle-Perron-Frobenius operator} $L_f$, defined by (in the case of the circle)
\begin{equation}\label{EqDefRPF}
L_f \phi : y \mapsto \sum_{f(x)=y} \frac{\phi(x)}{f'(x)}.
\end{equation}
The crucial remark is that if $\phi$ is the density of some measure $\mu$ on $\Sp^1$, then $L_f \phi$ is the density of the measure $f_*\mu$.
Ruelle proved that if $f$ is a $C^{1+\alpha}$ expanding map of the circle, then the functions $L_f^k(1)$ converge exponentially fast -- for the $C^\alpha$ topology -- towards the density of a measure called SRB (for Sinaï-Ruelle-Bowen). This measure is moreover the unique absolutely continuous $f$-invariant probability measure, and the unique physical measure\footnote{A measure $\mu$ is called \emph{physical} if there is a positive Lebesgue measure set of points $x\in\Sp^1$ whose Birkhoff sums $\frac{1}{n}\sum_{k=0}^{n-1}\delta_{f^k(x)}$ converge weakly towards $\mu$.} \cite{VianaStoch}. In summary, the measures $f^k_*(\Leb)$ will converge ``exponentially'' towards $SRB$. Together with Theorem~\ref{MainTheo}, this implies that there is some regime $k\ll N$ in which the measures $(f_N^k)_*(\Leb_N)$ converge towards $SRB$.

\begin{coro}
Let $r > 1$, and $f$ a generic $C^r$ expanding map of $\Sp^1$. Then there exists a constant $C>1$, depending only on $\inf f'$ and the $C^{r-1}$ norm of $f'$, such that for any $\varep>0$ and any $k\ge -C\log\varep$, there is $N_0\in\N$ such that for any $N\ge N_0$, one has
\[\Disc\Big(SRB, (f_N^k)_*(\Leb_N)\Big)\le\varep.\]
\end{coro}

This corollary can in fact be obtained from a direct computation without genericity assumption (see \cite[Theorem 12.19]{Guih-These} for an explicit statement). Theorem \ref{MainTheo} specifies this convergence by estimating the Cramér distance between continuous and discretized dynamics.
\bigskip

Of course, the discretization procedure studied in this paper is very primitive and inefficient to actually compute SRB measures for circle expanding maps, compared to some other algorithms like the Ulam approximation (see \cite{MR2863582} for a survey about this subject, and our second paper \cite{paper2}). Our aim here is to describe to what extent the very naive algorithm for computing SRB measures actually computes an approximation of this measure or not.

\subsection*{Overview of the paper}

The proof of Theorem \ref{MainTheo} follows the strategy of \cite{Gui15a}: as pointed out by Lanford in \cite{MR1678095}, it is very fruitful to mimic the proof strategy for the problem consisting in deducing the laws of fluid dynamics from a microscopic model of a gas (hydrodynamic limit), by introducing an intermediate mesoscopic space scale between the microscopic scale $1/N$ of the grid and the macroscopic scale. 

So we will start with the study of the corresponding local problem, i.e. the linear case (Sections \ref{SelSimplLinear} and \ref{Sectree}). First, we will	 define the corresponding notion of discretization for sequences of linear maps of $\R$. The first step is to link the Cramér distance with the roundoff errors made at each iteration (Proposition \ref{WasDistLin}). It turns out that under generic conditions, these roundoff errors are statistically uniformly distributed (Proposition \ref{Flocker}), a fact which had previously been obtained by Flockermann in his thesis \cite[Theorem~10 page 44]{Flocker}; we will here give two new proofs of this result in Subsections \ref{SecUnifEr1} and \ref{SecUnifEr2}. As in \cite{Gui15a}, both are based on ideas from the theory of quasicrystals; the first one follows a direct approach\footnote{This direct approach was already presented in the thesis \cite[Chapter 9]{Guih-These}.}, while the second one gives more precise results that allow to get a formula for the cumulated difference \eqref{DefDiscR} for time $k$ at some point $x\in\Z$ (which is used to define the Cramér distance). This formula involves the value of a piecewise linear functional on a $k$-dimensional torus at some point depending explicitly on $x$ (Proposition \ref{propDiscQuot}). This proposition is somehow the heart of the paper, as the first approach fails when one tries to pass to the more complicated framework of the tree linear case.

The set of time-$k$ preimages of a point $x\in\Sp^1$ by an expanding map $f$ has a structure of complete $d$-ary tree of height $k$, where $d$ is the degree of $f$. So the next step is to study a model of discretizations of linear maps that decorate such a complete $d$-ary tree (where the coefficients of the linear maps correspond to the derivatives of the expanding map). This is done in Section \ref{Sectree}, where we will deduce the behaviour of the Cramér distance in this framework from the study conducted in the previous section.

Finally, in Section \ref{SecExpand}, we will use these results to prove Theorem \ref{MainTheo}. It will be done in two steps. First, we will get a formula involving the derivatives of $f$ along paths of the preimage tree. It will be achieved by combining Thom's transversality theorem for generic maps with the study of the linear case and some suitable application of Taylor's formula. This third tool is elementary but rather technical, and will be obtained by applying a result of \cite{Gui15a}. The second step is quite elementary and will allow us to write the formula of Theorem \ref{ThWasDist} in the nicer way of Theorem \ref{MainTheo}.

\subsection*{Numerical experiments}

In \cite{paper2} we conduct some numerical experiments relative to the present article. Our aim is twofold: first, establish the time scale where Theorem~\ref{MainTheo} stays valid on some actual examples; second, for bigger times, try to determine numerically the phenomena underlying the behaviour of \eqref{EqStudied}.

It turns out that on the examples we tested, Theorem~\ref{MainTheo} stays valid until times $k$ typically logarithmic in $N$, and that in the regime where $k\ge \log N$ the evolution of \eqref{EqStudied} is not satisfyingly described by some model involving only random perturbations of the dynamics: the fact (involved in the prof of Theorem~\ref{MainTheo}) that orbits that merge stay together forever thereafter has a significant impact on the Cramér distance \eqref{EqStudied}. In \cite{paper2} we propose a model taking this phenomenon into account. We conjecture that this model captures sufficiently well the relevant features of discretizations in the middle-term range to approximate well the evolution of \eqref{EqStudied}.

\subsection*{Context}

To our knowledge, the first attempt of numerical approximation of physical measures of some Anosov diffeomorphisms dates back to le late 70's, with the works \cite{MR534103, MR0478237} where the authors study among others the Arnold cat map and some of its perturbations.

The idea of O.E. Lanford consisting in adjusting the length of orbit segments to the discretization order had already been developed in the late $80$'s by Abraham~Boyarsky. In \cite{MR862028}, he explains heuristically why one usually finds absolutely continuous measures on simulations. His arguments are based on the tracking of long segments of orbits; the only obstacle for the obtaining of a rigorous proof is the lack of uniformity in Birkhoff 's ergodic theorem\footnote{In fact, the intuition of Boyarsky works for a uniquely ergodic homeomorphism, as proved by Miernowski in \cite{MR2279269}. Unfortunately, his result is false in general (see \cite{MR3789171}).}.

In \cite{MR959419}, Boyarsky together with Pawel Góra establish the following result, which also relates to the obtaining of absolutely continuous measures from discretizations. \emph{Suppose that $f$ has a unique absolutely continuous invariant measure $\mu$, and that there exists $\alpha>0$ such that there is a subsequence of $f_N$ admitting a segment of orbit of length bigger than $\alpha \card(E_N)$. If we denote $\nu_N$ the uniform measure on this segment of orbit, then $\nu_N\rightharpoonup \mu$.}

The existence of an orbit segment of length proportional to that of the grid seems to be rarely verified (for example it is not true for a generic circle expanding map, simply because the degree of recurrence is zero). Despite this, this seems to be one of the first theoretical results about discretizations of dynamical systems.

The rigourous study of discretizations of \emph{generic} dynamics has first been proposed by \'Etienne Ghys in the large audience article \cite{Ghys-vari}. From this viewpoint, the case of circle homeomorphisms is now quite well understood, due to the work of Tomasz Miernowski \cite{MR2279269}, whose conclusion is essentially that the discretizations' dynamics resemble the homeomorphism's one (see also the recent preprint \cite{galatolo2021quantitative}). The higher dimensional case has been tackled by the first author in his thesis \cite{Guih-These}, which includes the case of generic homeomorphisms \cite{Guih-discr}, $C^1$ diffeomorphisms \cite{MR3789171} and $C^r$ diffeomorphisms and expanding maps\footnote{As already explained, the proofs of the present article are based on the strategy of this paper.} \cite{Gui15a}. In particular, the article \cite{MR3789171} exhibits the following quite unexpected phenomenon\footnote{A similar statement holds for generic measure-preserving $C^1$-diffeomorphisms.}.

\begin{theo*}[\cite{MR3789171}]
Take a point $x\in \Sp^1$, and a generic $C^1$ circle expanding map $f$. The orbit of $x$ under $f_N$ is finite thus eventually periodic; denote $\mu_N^x$ the uniform measure on the limit periodic orbit. Then the sequence $(\mu_N^x)_{N\ge 0}$ accumulates on the whole set of $f$-invariant probability measures.
\end{theo*}

Hence, the ergodic behaviour under the discretization $f_N$ of a (Baire) typical point of the circle does not converge towards the unique physical measure (which exists and is singular, see \cite{MR1845327,MR1688216}). However, one can expect some more convergent behaviour by averaging over $x\in \Sp^1$. This leads to the following question.

\begin{ques}
For $r>1$ and a generic $C^r$  expanding map of $\Sp^1$, do the measures\footnote{The convergence of these measures in $k$ is ensured by the finiteness of the map $f_N$.}
\[\lim_{k\to +\infty} \frac{1}{k} \sum_{i=0}^{k-1} (f_N^i)_*(\Leb_N)\]
converge, when the parameter $N$ goes to infinity, towards the SRB measure of $f$?
\end{ques}

This question seems out of reach using only the techniques used in this paper. In the second article of this series \cite{paper2}, we show numerical simulations that suggest that the answer to this question may be yes in general.

\subsection*{Acknowledgements}

This project was partially supported by a PEPS/CNRS project and the ANR CODYS. The authors warmly thank Nina Heloin for his careful reading of a first version of this text, and Djalil Chafai for the references about the name of our distance $\disc$.

\section{Preliminaries: distances on measures}\label{SecDistMeas}

Let $\mu$ and $\nu$ be two probability measures defined on $\Sp^1 = \R/\Z$, identified with $[0,1[$. Let $F$ and $G$ be their respective cumulative distribution functions, and $H=F-G$

\begin{definition}
The $L^1$ Wasserstein distance between $\mu$ and $\nu$ can be defined by the formula\footnote{In \cite{MR1335789} it is explained why this formula coincides with the classical definition of the Wasserstein distance.}
\[W_1(\mu,\nu) = \min_{c\in\R} \int_0^1 |H-c|,\]
and the minimum is realized by the median of $F-G$, i.e. 
\[c_0 = \frac12\Big(\sup\big\{c\mid \Leb(H<c) < \Leb(H>c)\big\} + \inf\big\{c\mid \Leb(H<c) > \Leb(H>c)\big\} \Big).\]

Similarly we define another distance on the set of probability measures on $\Sp^1$, which we call \emph{Cramér distance}:
\begin{equation}\label{DefDisc}
\Disc(\mu,\nu) = \left(\min_{c\in\R} \int_0^1 (H(x)-c)^2\ud x\right)^{1/2},
\end{equation}
and the minimum is realized by the mean of $F-G$, i.e. by the number $c_1 = \int_0^1 H$:
\[\Disc(\mu,\nu) = \left(\int_0^1 \left(H(x) - \Big(\int_0^1 H(y)\ud y\Big)\right)^2 \ud x\right)^{1/2}.\]
\end{definition}

Remark that this last expression looks like a variance, and this fact will be useful in the sequel. The proofs of the statements about the numbers $c$ realizing the minima are simple and left to the reader.

\begin{lemme}\label{LemSpanWeakStar}
The Cramér distance $\Disc$ is a distance which is invariant under translation. Moreover,
\[W_1(\mu,\nu) \le \Disc(\mu,\nu)\]
and
\[\Disc(\mu,\nu) \le \sqrt 2 W_1(\mu,\nu)^{1/2}.\]
Thus $W_1$ and $\Disc$ span the same topology (the weak-*).
\end{lemme}

\begin{proof}[Proof of Lemma \ref{LemSpanWeakStar}]
First we prove that $\Disc$ is invariant under translation. More precisely, for any $a\in [0,1]$, we let $F_a$ and $G_a$ be the cumulative distribution functions of $\mu$ and $\nu$ seen as measures on $[a,a+1]$, and set $H_a = F_a-G_a : [a,a+1] \to \R$. What we want to prove is that for any $a\in[0,1]$, one has 
\[\Disc(\mu,\nu)^2 = \int_a^{a+1} \left(H_a(x) - \Big(\int_a^{a+1} H_a\Big)\right)^2 \ud x.\]
One has $H_a(u) = \int_a^u \ud(\mu-\nu)$, thus, using the fact that $\int_0^1 \ud(\mu-\nu) = 0$,
\[H_a(u) = \left\{\begin{array}{ll}
\int_0^u \ud(\mu-\nu) - \int_0^a \ud(\mu-\nu) = H(u)-H(a) & \text{if } u\le 1\\
\int_1^u\ud(\mu-\nu) + \int_a^1 \ud(\mu-\nu) = H(u-1)-H(a) & \text{if } u\ge 1.
\end{array}\right.\]
Hence, $\int_a^{a+1} H_a = \int_0^1 (H - H(a))$, which implies that 
\[H_a(u) - \int_a^{a+1} H_a = H(u \mod 1) - \int_0^1 H,\]
and thus
\[\int_a^{a+1} \left(H_a(x) - \Big(\int_a^{a+1} H_a\Big)\right)^2 \ud x = \int_0^1 \left(H(x) - \Big(\int_0^1 H\Big)\right)^2 \ud x.\]
\medskip

We now come to the proof of the inequalities. The first one is simply Cauchy-Schwarz inequality applied to the map $H-c_1$.

For the second one, remark that $H\in [-1,1]$, so that $c_0 \in [-1,1]$ and $|H-c_0|\le 2$. Hence, $(H-c_0)^2 \le 2|H-c_0|$ and 
\[\Disc(\mu,\nu)^2 \le \int_0^1 (H-c_0)^2 \le 2\int_0^1 |H-c_0| = 2 W_1(\mu,\nu).\]
\end{proof}


\section{The simple linear case}\label{SelSimplLinear}

Recall that the goal of this section is to treat the corresponding problem of Cramér distance between actual and discretized systems, for sequences of linear maps. We first set definitions, which are made to mimic the ones on the circle for the case of linear maps of $\R$.

\subsection{Definitions}

We denote $\N = \{0,1,2,\dots\}$.

\subsubsection*{Discretizations of linear maps.}
Let $(\ell_m)_{m\ge 1}$ be a sequence of homotheties of $\R$, of parameters $\lambda_m>1$, i.e. $\forall y\in\R$, $\ell_m(y) = \lambda_m y$. We fix $k\in\N^*$.

\begin{definition}
The \emph{discretization} of a linear map $\ell : \R\to \R$ is the map $\widehat \ell : \Z\to\Z$ such that for any $x\in\Z$, $\widehat \ell (x)$ is the integer closest to $\ell(x)$. More precisely, $\widehat \ell (x)$ is the unique integer such that 
\[\widehat \ell (x) - \ell(x) \in \left]-\frac12, \frac12\right].\]
We will denote
\begin{equation}
\Ell^k =  \ell_{k} \circ \cdots \circ \ell_1, \qquad \widehat{\Ell}^k = \widehat \ell_{k} \circ \cdots \circ \widehat \ell_1 \qquad \text{and} \qquad \widetilde\lambda_m = \prod_{i=m+1}^k \lambda_i
\end{equation}
(with the convention that $\widetilde\lambda_{k} = 1$).
\end{definition}

\subsubsection*{Expectation, covariance.}
The expectation $\mathbb E$ and the covariance $\operatorname{Var}$ of a map $\mathcal E : \N\to\R$ are defined by (whenever the limits make sense)
\begin{equation}\label{EqVariance}
\mathbb E [\Er] =  \lim_{R\to +\infty} \frac{1}{R}\sum_{x=0}^{R-1} \mathcal E_x\qquad \text{and} \qquad
\Var (\Er) =  \lim_{R\to +\infty} \frac{1}{R}\sum_{x=0}^{R-1} \big(\mathcal E_x-\mathbb E[\mathcal E]\big)^2.
\end{equation}

\subsubsection*{Cumulated difference and Cramér distance on $\R$.}
Let $\mu$ and $\nu$ be two measures on $\R$. Their \emph{cumulated difference} at $y>0$ is the number 
\begin{equation}\label{DefDiscR}
c\delta_y(\mu,\nu) = \mu\big(]0,y]\big) + \frac12\mu\big(\{0\}\big) - \nu\big(]0,y]\big) - \frac12\nu\big(\{0\}\big).
\end{equation}

We define the \emph{($L^2$-) Cramér distance} $\Disc$ as the $L^2$-average of the cumulated difference $c\delta$ \eqref{DefDiscR} (when the limit exists):
\[\Disc(\mu,\nu) = \lim_{R\to +\infty}\Disc_R(\mu,\nu)\qquad
\text{where}
\qquad \Disc_R(\mu,\nu) =  \left(\frac{1}{R}\int_0^R c\delta_y(\mu,\nu)^2\right)^{1/2}.\]

We will be interested in the case where the measures $\mu$ and $\nu$ are respectively:
\begin{itemize}
\item the (correctly normalized) Lebesgue measure $\widetilde\lambda_0^{-1}\Leb $;
\item the uniform measure on the image set $\widehat \Ell^k(\Z)$, that is $\sum_{n\in\Z} \delta_{\widehat \Ell^k(n)}$
\end{itemize}
In the sequel we will denote, when no confusion is possible,
\begin{align*}
c\delta_y & \overset{\text{def.}}{=} c\delta_y\left(\widetilde\lambda_0^{-1}\Leb,\ \sum_{n\in\Z} \delta_{\widehat \Ell^k(n)} \right) \nonumber \\
         & = \frac{y}{\widetilde\lambda_0} - \card\big\{x\in\N\mid \widehat{\Ell}^k (x)\le y\big\} + \frac{1}{2},
\end{align*}
and the same for the discrepancies $\Disc_R$ and $\Disc$.

The half weight given to the singleton $\{0\}$ restores symmetry and ensures that the map $c\delta$ has zero mean (see Remark \ref{RemNormaliz}).

\subsection{Roundoff errors}

%

The roundoff error made at the $m$-th iteration is defined as the difference between the images of $\Ell^{m-1}(x)$ by the discretization $\widehat \ell_m$ and the initial map $\ell_m$, that is
\[e_x^m = \big(\widehat \ell_m - \ell_m\big)\big(\widehat\Ell^{m-1}(x)\big) \in ]-1/2,1/2].\]

The distribution of the vectors
\[\Varep_x^k \overset{\text{def.}}{=} (e_x^1,\dots,e_x^k)\]
when $x$ ranges over $\Z$ is given by the following proposition  due to P.~P. Flockermann (see the thesis \cite{Flocker}, Theorem~10 page 44). We will give two alternative proofs of this proposition, both based on linear algebra (contrary to the original proof of Flockermann).

In the following, by \emph{$\Q$-free family} we mean a family of real numbers such that, when adding 1 to this family, one gets a free family when $\R$ is seen as a $\Q$-vector space. A finite family of real numbers $(\lambda_i)_{1\le i\le k}$ is $\Q$-free if and only if there is no nonzero family $(m_i)_{1\le i\le k}$ of integers satisfying $\sum_{i=1}^k m_i\lambda_i\in\Z$.

\begin{prop}[Flockermann]\label{Flocker}
If the family $(\widetilde\lambda_m^{-1})_{0\le m\le k}$ is $\Q$-free (which is a generic condition on the $\lambda_i$'s), then the roundoff error vectors $(\Varep_x^k)_{x\in\Z}$ are equidistributed in $]-1/2,1/2]^k$.
\end{prop}

We will give two independant proofs of this proposition in Sections \ref{SecUnifEr1} and \ref{SecUnifEr2}: the first one will be a first step to understand the general proof strategy, the second one will be a byproduct of the proof of Proposition~\ref{propDiscQuot}.

From the roundoff errors $\Varep_x^k$ it is possible to deduce the global error
\[\Er_x^k = \widehat \Ell^k(x) - \Ell^k(x)\]
made after $k$ iterations. Indeed, we have 
\begin{align*}
\Er_x^{k+1} = & \big(\widehat \ell_{k+1} - \ell_{k+1}\big)\big(\widehat\Ell^k(x)\big) + \ell_{k+1} \big( \widehat\Ell^k(x) - \Ell^k(x)\big) \\
= & e_x^{k+1} + \ell_{k+1}(\Er_x^k).
\end{align*}
From this recurrence relation, we deduce that
\begin{equation}\label{FormulEx}
\Er_x^k = \sum_{m=1}^k \widetilde\lambda_{m} e_x^m.
\end{equation}
Recall that Proposition \ref{Flocker} ensures that when the family $( \widetilde\lambda_m^{-1})_{1\le m\le k}$ is $\Q$-free, then the errors $e_x^m$ are independent and identically distributed in $[-1/2,1/2]$ (because $\Varep_x^k$ is equidistributed on the product space $[-1/2, 1/2]^k$). From that we deduce the law of the global error $\Er^k$, and in particular its covariance
\begin{equation}\label{FormulEx2}
\operatorname{Var} (\Er^k) = \sum_{m=1}^k \widetilde\lambda_{m}^2 \operatorname{Var}(e^m) = \frac{1}{12}\sum_{m=1}^k \widetilde\lambda_{m}^2.
\end{equation}
In particular, if there exists $\alpha>1$ such that $\lambda_m\ge \alpha$ for every $m$, then 
\[\operatorname{Var} (\mathcal E^k) \ge \frac{1}{12}\sum_{m=1}^k \alpha^{2m} = \frac{\alpha^2 (\alpha^{2k}-1)}{12(\alpha^2-1)}.\]

\subsection{Cramér distance and roundoff errors}

In this subsection, we link the asymptotic behaviours of the Cramér distance $\Disc$ with that of the roundoff errors. We state all properties before proving them. Note that some of these proofs will use Proposition \ref{Flocker}, which will be proved later on without using the results of this subsection.

The first lemma says that the mean of the map $c\delta$ is zero.

\begin{lemme}\label{MeanZero}
\[\mathbb E[c\delta_{x+1/2}] = \lim_{R\to +\infty}\frac{1}{R}\int_0^R c\delta_y \ud y = 0.\]
\end{lemme}

Remark that we have two different notions of means: one continuous (with an integral) and one discrete (with a sum). Both can be easily related: the following lemma links the Cramér distance $\Disc$ which is obtained as a continuous average of the cumulated difference $c\delta$, with the variance of the map $c\delta$ taken on half integers. 

\begin{lemme}\label{LemCalcDiscdisc}
Whenever the Cramér distance and the variance make sense, 
\begin{equation}\label{EqCalcDiscdisc}
\Disc^2 = \frac{1}{12 \widetilde\lambda_0^2} + \Var_{x\in\Z} \left(c\delta_{x+\frac12}\right).
\end{equation}
\end{lemme}

Finally, the following proposition deals with the covariance: it links the average Cramér distance $\Disc$ with the covariance of $\Er$.

\begin{prop}\label{WasDistLin}
Let $k\in\N$, and a family $(\ell_m)_{1\le m\le k}$ of 	 of $\R$ of parameters $(\lambda_m)_{1\le m \le k}$ strictly bigger than 1. Whenever the Cramér distance and the variance make sense, 
\[\Disc^2  = \frac{1}{12} + \frac{1}{\widetilde\lambda_0^2} \operatorname{Var}(\Er^k).\]
\end{prop}

Note that the factor $1/12$ corresponds to the covariance of the uniform distribution on the interval $[-1/2,1/2]$.

Combined with Proposition \ref{Flocker} (more precisely, Equation \eqref{FormulEx2}), this immediately gives the following corollary.

\begin{coro}\label{CoroWasDistLin}
If the family $(\widetilde\lambda_m)_{0\le m\le k}$ is $\Q$-free, then
\[\Disc^2 = \frac{1}{12\widetilde\lambda_0^2} \sum_{m=0}^{k}\widetilde\lambda_m^2.\]
\end{coro}

\begin{proof}[Proof of Lemma~\ref{MeanZero}]
The first equality comes from the fact that the map $c\delta$ is affine with slope $\widetilde\lambda_0^{-1}$ in restriction to any interval which contains no integer. We will see a more detailed proof of a very similar fact during the proof of Lemma \ref{LemCalcDiscdisc}.

We are left to prove the second equality. Let
\[\Er'_x = \widetilde\lambda_0^{-1} \Er_x^k.\]
Fix $R>0$, set $R'=\widetilde\lambda_0 R$, and denote $x_0 = x_0(R)$ the biggest integer $x$ such that $x+\Er'_x \le R$. A linear change of variables leads to
\begin{align*}
\frac{1}{R'}\int_0^{R'} c\delta_{y'} \ud y' 
& = \frac{1}{R'}\int_0^{R'}\left( \frac{y'}{\widetilde{\lambda}_0} + \frac12 - \sum_{x=0}^{x_0}\1_{\ell^k(x) + \Er_x\le y'}\right)\ud y' \\
& = \frac{1}{R} \int_0^R \left(y+\frac12 - \sum_{x=0}^{x_0} \1_{x+\Er'_x \le y}\right)\ud y.
\end{align*}
Hence, 
\begin{align*}
\frac{1}{R'}\int_0^{R'} c\delta_{y'} \ud y'
& = \frac{1}{R}\left( \frac{R^2}{2} + \frac{R}{2} - \sum_{x=0}^{x_0} \int_0^R \1_{y\ge x+\Er'_x} \ud y\right)\\
& =  \frac{R}{2} + \frac{1}{2} - \frac{1}{R}\sum_{x=0}^{x_0} (R-x-\Er'_x)\\
& =  \frac{R}{2} + \frac{1}{2} - (x_0+1) + \frac{x_0(x_0+1)}{2R} + \frac{1}{R}\sum_{x=0}^{x_0} \Er'_x\\
& =  \frac{(R-x_0)^2 - (R-x_0)}{2R} + \frac{x_0}{R}\frac{1}{x_0}\sum_{x=0}^{x_0} \Er'_x.
\end{align*}
But $|R-x_0|$ is uniformly bounded on $R$ (because $\Er'_x$ is uniformly bounded on $R$), so the first term tends to 0 when $R'$ goes to infinity. The second term, for itself, tends to the mean of $x\mapsto \Er'_x$, which is zero by Proposition \ref{Flocker}.
\end{proof}

\begin{proof}[Proof of Lemma \ref{LemCalcDiscdisc}]
The proof simply consists in remarking that the map $c\delta$ is affine with slope $\widetilde\lambda_0^{-1}$ in restriction to any interval that contains no integer. When $R\in\N$, one has
\begin{align*}
\frac{1}{R} \int_0^R c\delta_y^2 \ud y
= & \frac{1}{R}\sum_{x=0}^{R-1}\int_{-1/2}^{1/2} \left(c\delta_{x+\frac12} + \frac{y}{\widetilde\lambda_0}\right)^2 \ud y\\
= & \frac{1}{R}\sum_{x=0}^{R-1} c\delta_{x+\frac12}^2 + \frac{1}{12 \widetilde\lambda_0^2}.
\end{align*}
But by Proposition \ref{MeanZero}], one has $\mathbb E\big[c\delta_{k+\frac12}\big] = 0$; this gives directly the lemma.
\end{proof}

\begin{proof}[Proof of Proposition~\ref{WasDistLin}]
We reuse the notations of proof of Lemma \ref{MeanZero}: we set $\Er'_x = \widetilde\lambda_0^{-1} \Er_x$, $R'=\widetilde\lambda_0 R$, and denote $x_0$ the biggest integer $x$ such that $x+\Er'_x \le R$.

We will prove that
\[\Disc_R^2 \left(\widehat{\Ell}^k(\Z)\, , \, \widetilde\lambda_0^{-1} \Leb\right)\underset{R\to +\infty}{\longrightarrow} \frac{1}{12} + \operatorname{Var}(\Er').\]

By Equation \eqref{DefDiscR}, we have
\[\Disc_{R'}^2 = \frac{1}{R'} \int_0^{R'} \left(y'\widetilde\lambda_0^{-1} - \sum_{x=0}^{x_0} \1_{\widetilde\lambda_0 x+\Er_x \le y'} + \frac12\right)^2\ud y'.\]
A linear change of variables leads to
\begin{align*} 
\Disc_{R'}^2
		& = \frac{1}{R} \int_0^R \left(y+\frac12 - \sum_{x=0}^{x_0} \1_{x+\Er'_x \le y}\right)^2\ud y\\
		& = \frac{1}{R} \int_0^R \left(\left(y+\frac12\right)^2 - 2\left(y+\frac12\right)\sum_{x=0}^{x_0} \1_{x+\Er'_x \le y} + \sum_{x,x'=0}^{x_0} \1_{x+\Er'_x \le y} \1_{x'+\Er'_{x'} \le y}\right)\ud y\\
		& = \frac{1}{R} \Bigg[\int_0^R \left(y+\frac12\right)^2 \ud y - 2\sum_{x=0}^{x_0} \int_0^R \left(y+\frac12\right) \1_{x+\Er'_x \le y}\ud y \\
		& \phantom{= \frac{1}{R} \Bigg[}+ \sum_{x,x'=0}^{x_0} \int_0^R \1_{x+\Er'_x \le y} \1_{x'+\Er'_{x'} \le y}\ud y\Bigg]\\
		& = \frac{1}{R} \left[\int_0^R \left(y+\frac12\right)^2 \ud y - 2\sum_{x=0}^{x_0} \int_{x+\Er'_x}^R \left(y+\frac12\right) \ud y + \sum_{x,x'=0}^{x_0} \int_{\max\big(x+\Er'_x, x'+\Er'_{x'} \big)}^R 1\ud y\right]\\
		& = \frac{1}{R} \Bigg[\frac{\left(R+\frac12\right)^3}{3} - \frac{\left(\frac12\right)^3}{3} - \sum_{x=0}^{x_0} \left(\left(R+\frac12\right)^2 - \left(x+\Er'_x+\frac12\right)^2\right)\\
		& \phantom{= \frac{1}{R} \Bigg[} + \sum_{x,x'=0}^{x_0}\Big( R-\max\big(x+\Er'_x,x'+\Er'_{x'}\big)\Big) \Bigg].
\end{align*}
But by construction, the map $x\mapsto x+\Er'_x$ is increasing, so the last sum can be reindexed:
\begin{align*}
\sum_{x,x'=0}^{x_0}\Big(R-\max\big(x+\Er'_x,x'+\Er'_{x'}\big)\Big) & = \sum_{m=0}^{x_0}\ \sum_{\substack{\max(x,x')=m \\ x,x'\ge 0}} \big(R-m-\Er'_m\big)\\
& = \sum_{m=0}^{x_0} (2m+1)\big(R-m-\Er'_m\big),
\end{align*}
so one gets:
\begin{align*}
\Disc_{R'}^2 = & \frac{1}{R} \Bigg[\frac{R^3}{3} +\frac{R^2}{2} + \frac{R}{4} - \sum_{x=0}^{x_0} \bigg(R^2 + R + \frac14 - \Big(x^2 + {\Er'_x}^2 + \frac14 + 2x\Er'_x + x + \Er'_x \Big)\bigg)\\
          & + \sum_{x=0}^{x_0} (2x+1)\big( R-(x+\Er'_x)\big)\Bigg]\\
        = & \frac{1}{R} \Bigg[\frac{R^3}{3} +\frac{R^2}{2} + \frac{R}{4} + \sum_{x=0}^{x_0} \bigg(-R^2-x^2 + 2xR + {\Er'_x}^2 \bigg)\Bigg]\\
        = & \frac{1}{R} \Bigg[\frac{R^3}{3} +\frac{R^2}{2} + \frac{R}{4} -R^2(x_0+1) - \left(\frac{x_0^3}{3} + \frac{x_0^2}{2} + \frac{x_0}{6}\right) + 2R\frac{x_0(x_0+1)}{2} + \sum_{x=0}^{x_0} {\Er'_x}^2\Bigg]\\
	    = & \frac{1}{R} \left[ \frac{(R-x_0)^3}{3} - \frac{(R-x_0)^2}{2} + \frac{R-x_0}{6} \right] + \frac{1}{12} + \frac{1}{R} \sum_{x=0}^{x_0} {\Er'_x}^2.
\end{align*}
As in the proof of Lemma \ref{MeanZero}, we have that $|R-x_0|$ is uniformly bounded in $R$, so that the first term tends to $0$ as $R'$ goes to infinity. As the mean of $x\mapsto \Er'_x$ is 0 (see Proposition \ref{RoundoffLin}), we get finally that 
\[\lim_{R'\to +\infty} \Disc_{R'}^2 = \frac{1}{12} + \Var(\Er').\]
\end{proof}

\subsection{Uniform distribution of errors: first proof}\label{SecUnifEr1}

This subsection presents the first proof of uniform distribution of errors. It is easier than the one we will see in the next section, and consists in computing projections on $k$-dimensional tori of vectors depending on the initial condition $x$.

Fix $k\ge 0$. Recall that we have $\Varep_x^k = \big(e_x^1,\cdots,e_x^k\big)$. Moreover, we set $\widehat\Ell_x = \big(\widehat\Ell^1(x),\cdots,\widehat\Ell^k(x)\big)$ the vectors made of the $k$ firsts iterates of $x$ under the discretizations, and denote $u_x = (\lambda_1x,0^{k-1})\in\R^{k}$, $u_\Z = (\lambda_1\Z,0^{k-1})\subset\R^{k}$, $W^k = ]-1/2,1/2]^{k}$ and
\[N_{\lambda_1,\cdots,\lambda_k} = \begin{pmatrix}
-1  &     &        &     & \\
\lambda_2 & -1  &        &     & \\
    & \lambda_3 & \ddots &     & \\
    &     & \ddots & -1  & \\
    &     &        & \lambda_k & -1
\end{pmatrix}\in M_{k}(\R).\]
Finally, we denote $\pr_{W^k}$ the projection from $\R^k$ onto the fundamental domain $W^k$ of the quotient space $\R^k/N_{\lambda_1,\cdots,\lambda_k}\Z^k$ (this is indeed a fundamental domain because the matrix $N_{\lambda_1,\cdots,\lambda_k}$ is lower triangular with $-1$ on the diagonal; remark that in this case the matrix satisfies the conclusion of Haj\'os theorem \cite{MR0006425}).

The following proposition expresses the roundoff error vector $\Varep_x^k$ in terms of the projection of the vector $u_x$ on the fundamental domain $W^k$ of $\R^k/N_{\ell_1,\cdots,\ell_k}\Z^k$.

\begin{prop}\label{RoundoffLin}
\[\Varep_x^k = \pr_{W^k}(u_x).\]
Thus, when $x$ ranges over $\Z$, the roundoff error vectors $\Varep_x^k$ equidistribute on the set $\overline{\pr_{W^k}(u_\Z)}$.

In particular, as this set is symmetric with respect to 0, the means of each function $x\mapsto e_x^m$, and of $x\mapsto \Er_x$, is zero.
\end{prop}

\begin{proof}[Proof of Proposition \ref{RoundoffLin}]
As 
\[N_{\lambda_1,\cdots,\lambda_k} \widehat\Ell_x = \begin{pmatrix}
-\widehat\Ell^1(x) \\
\lambda_2 \widehat\Ell^1(x) - \widehat\Ell^2(x)\\
\lambda_3 \widehat\Ell^2(x) - \widehat\Ell^3(x)\\
\vdots\\
\lambda_k \widehat\Ell^{k-1}(x) - \widehat\Ell^{k}(x)
\end{pmatrix}
= \begin{pmatrix}
\lambda_1 x - e_x^1\\
-e_x^2\\
-e_x^3\\
\vdots\\
-e_x^k
\end{pmatrix} = u_x - \Varep_x^k,\]
the vector $u_x$ can be decomposed into
\begin{equation}\label{decompfund}
u_x  = N_{\lambda_1,\cdots,\lambda_k}\widehat\Ell_x + \Varep_x^k,
\end{equation}
with $\widehat\Ell_x \in \Z^{k}$ and $\Varep_x^k \in W^k$ (recall that $W^k = ]-1/2,1/2]^{k}$). As $W^k$ is a fundamental domain of $N_{\lambda_1,\cdots,\lambda_k}\Z^{k}$, this is a decomposition of $u_x$ into the sum of an element of the lattice $N_{\lambda_1,\cdots,\lambda_k}\Z^{k}$ and an element of a fundamental domain of this lattice.

The vector $u_x$ being fixed, this condition characterizes completely $\Varep_x^k$ and $\widehat\Ell_x$. In particular, $\Varep_x^k$ is equal to the projection of $u_x$ on $W^k$ modulo $N_{\lambda_1,\cdots,\lambda_k}\Z^{k}$. This implies that the roundoff error vectors $\Varep_x^k$ equidistribute on the set $\overline{\pr_{W^k}(u_\Z)}$ when $x$ ranges over $\Z$.
\end{proof}

Let us explain how this proposition implies Proposition~\ref{Flocker}.

\begin{proof}[Proof of Proposition~\ref{Flocker}]
We begin by remarking that by \eqref{decompfund}, $N_{\lambda_1,\cdots,\lambda_k}^{-1}\Varep_x^k$ is equal to the projection of $N_{\lambda_1,\cdots,\lambda_k}^{-1} u_x$ on $N_{\lambda_1,\cdots,\lambda_k}^{-1} W^k$ modulo $\Z^{k}$ (remark that $N_{\lambda_1,\cdots,\lambda_k}^{-1} W^k$ is a fundamental domain of $\Z^k$). This implies that the sequence of errors $\Varep_x^k$ is equidistributed in $\R^k/\Z^k$ if and only if the vectors $N_{\lambda_1,\cdots,\lambda_k}^{-1} u_x$ are equidistributed modulo $\Z^{k}$ when $x$ ranges over $\Z$. For this purpose, the matrix $N_{\lambda_1,\cdots,\lambda_k}^{-1}$ can be easily computed:
\[N_{\lambda_1,\cdots,\lambda_k}^{-1} = \begin{pmatrix}
-1  &      &        &        & \\
-\lambda_2 & -1 &    &        & \\
-\lambda_3\lambda_2  & -\lambda_3 & \ddots &        & \\
\vdots   & \vdots & \ddots & -1 & \\
{}\ -\lambda_k \cdots \lambda_2\ {} & {}\ -\lambda_k \cdots \lambda_3\ {} & \cdots & -\lambda_k & -1
\end{pmatrix},\]
thus
\[
N_{\lambda_1,\cdots,\lambda_k}^{-1} u_x = -\begin{pmatrix}
\lambda_1\\
\lambda_2 \lambda_1\\
\lambda_3 \lambda_2 \lambda_1\\
\vdots\\
\lambda_k \cdots \lambda_1
\end{pmatrix}x
= - \widetilde \lambda_0\begin{pmatrix}
\widetilde\lambda_1^{-1}\\
\widetilde\lambda_2^{-1}\\
\vdots\\
\widetilde\lambda_{k}^{-1}
\end{pmatrix}x.
\]
As a consequence, by Weyl's criterion, the sequences of errors $\Varep_x^k$ is equidistributed in $W^k$ if and only if the family $( \widetilde\lambda_m^{-1})_{0\le m\le k}$ is $\Q$-free.
\end{proof}

\subsection{Cramér distance: a direct approach}\label{SecUnifEr2}

In the last subsection, we were given $x\in\Z$ and computed the roundoff errors along the positive orbit of $x$. Now, we adopt a different viewpoint: we are given $n\in \Z$ and want to determine whether $n$ belongs to $\widehat{\Ell}^k(\Z)$ or not; in the latter case we also want to determine the sequence of roundoff errors in the backward orbit of $n$. As in the previous section, we will see that these quantities only depend on the projection on some torus of a vector depending only on $n$. This will allow to compute the cumulated difference $c\delta_n$ from this projection.
\bigskip

\subsubsection*{Notations}
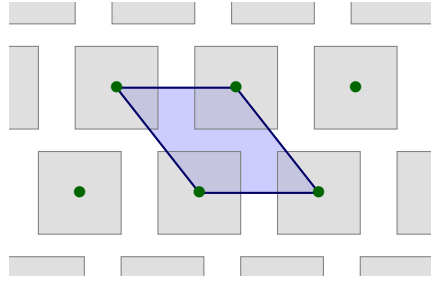
\begin{figure}[t]
\begin{center}
\begin{tikzpicture}[scale=1.1]
\clip (-2.3,-1) rectangle (2.8,2.3);
\fill[color=blue!20!white] (0,0) -- (1.445,0) -- (0.445,1.272) -- (-1,1.272) -- cycle;
\foreach\i in {-2,...,3}{
\foreach\j in {-1,...,3}{
\fill[color=gray!50!white,opacity = .5] (1.445*\i - \j-.5,1.272*\j-.5) rectangle (1.445*\i - \j+.5,1.272*\j+.5);
\draw[color=gray] (1.445*\i - \j-.5,1.272*\j-.5) rectangle (1.445*\i - \j+.5,1.272*\j+.5);
}}
\draw[color=blue!40!black,thick] (0,0) -- (1.445,0) -- (0.445,1.272) -- (-1,1.272) -- cycle;
\foreach\i in {-2,...,3}{
\foreach\j in {-1,...,3}{
\draw[color=green!40!black] (1.445*\i - \j,1.272*\j) node {$\bullet$};
}}
\end{tikzpicture}
\caption{The green dots are the points of $\widetilde \Lambda_k = \widetilde M_{\lambda_1,\cdots,\lambda_k}\Z^k$, the gray squares are $W^k + \widetilde \Lambda_k$ and the blue parallelogram is a fundamental domain of $X_k = \R^k/\widetilde\Lambda_k$.}\label{DessinRecap}
\end{center}
\end{figure}

Recall that $W^k = ]-1/2,1/2]^k$. We denote $\Lambda_k = M_{\lambda_1,\cdots,\lambda_k}\Z^{k+1}$, with
\[
M_{\lambda_1,\cdots,\lambda_k} = \left(\begin{array}{ccccc}
\lambda_1 & -1 &        &        & \\
    & \lambda_2  & -1   &        & \\
    &      & \ddots & \ddots & \\
    &      &        & \lambda_k    & -1\\
    &      &        &        & 1
\end{array}\right)\in M_{k+1}(\R),
\]
and $\widetilde \Lambda_k = \widetilde M_{\lambda_1,\cdots,\lambda_k}\Z^k$, with
\begin{equation}\label{DefMat}
\widetilde M_{\lambda_1,\cdots,\lambda_k} = \left(\begin{array}{ccccc}
\lambda_1 & -1 &        &         & \\
    & \lambda_2  & -1   &         & \\
    &      & \ddots & \ddots  & \\
    &      &        & \lambda_{k-1} & -1\\
    &     &        &          & \lambda_k
\end{array}\right)\in M_{k}(\R).
\end{equation}
(see Figure \ref{DessinRecap}). Finally, we denote $X_k = \R^k/\widetilde\Lambda_k$ the quotient space and $\pr_{X_k}$ the projection from $\R^k$ onto $X_k$. Remark that $X_k$ is a $k$-dimensional flat torus.

\bigskip

We begin by giving an alternative construction of the image sets $\widehat{\Ell}^k(\Z)$ in terms of model sets (see \cite{Gui15a}). Indeed, denoting $p_1$ the projection on the $k$ first coordinates and $p_2$ the projection on the last coordinate
\begin{align}
\widehat{\Ell}^k (\Z)\nonumber 
     & = \big\{p_2(\lambda)\mid \lambda\in \Lambda_k,\, p_1(\lambda)\in W^k\big\}\nonumber\\
     & = p_2\Big(\Lambda_k \cap \big(p_1^{-1}(W^k)\big)\Big).\label{CalcGamma}
\end{align}

Let us explain this construction. Fix a linear map $\ell : \R\to\R$ associated to $\lambda>1$. An integer $y\in\Z$ belongs to $\widehat\ell(\Z)$ iff there exists $x\in\Z$ such that $|\ell(x) - y|\le 1/2$. The last condition can be rephrased as $p_1(v) \in [-1/2,1/2]$, with 
\[v = \begin{pmatrix} \lambda & -1 \\ 0 & 1 \end{pmatrix} \begin{pmatrix} x \\ y \end{pmatrix}.\]
For two linear maps $\ell_1,\ell_2 : \R\to\R$ associated to $\lambda_1,\lambda_2>1$, a number $y\in\Z$ belongs to $(\widehat\ell_2\circ\widehat\ell_1)(\Z)$ iff there exist $x_1,x_2\in\Z$ such that $|\ell_2(x_2) - y|\le 1/2$ and $|\ell_1(x_1) - x_2|\le 1/2$ (and in this case $(\widehat\ell_2\circ\widehat\ell_1)(x_1) = \widehat\ell_2(x_2) = y$). These conditions can be rephrased as $p_1(v) \in [-1/2,1/2]^2$, with 
\[v = \begin{pmatrix} \lambda_1 & -1 & 0 \\ 0 & \lambda_2 & -1 \\ 0 & 0 & 1 \end{pmatrix} \begin{pmatrix} x_1 \\ x_2 \\ y \end{pmatrix}.\]
remark that in this case, the roundoff error is given by
\[\Varep_x^2 = \begin{pmatrix} x_2-\lambda_1 x_1 \\ y-\lambda_2 x_2\end{pmatrix} = -p_1(v),\]
and that $p_2(v) = y$. The same reasoning in arbitrary time $k$ leads to Equation \eqref{CalcGamma}.
\bigskip

Taking advantage from this viewpoint, we get the following proposition.

\begin{prop}\label{RoundoffLin2}
\[y\in\widehat{\Ell}^k (\Z) \quad \iff \quad y\in\Z\ \text{and}\ \pr_{X_k}(0^{k-1},y)\in \pr_{X_k}(- W^k).\]

In this case, if we denote by $w\in W^k$ the unique point satisfying $\pr_{X_k}(0^{k-1},y) = \pr_{X_k}(w)$, and $x\in\Z$ the unique integer such that $\widehat{\Ell}^k(x) = y$, then the roundoff errors satisfy $\Varep_x^k = -w$. As a corollary we get Proposition \ref{Flocker}.
\end{prop}

\begin{proof}[Proof of Proposition \ref{RoundoffLin2}]
By Equation~\eqref{CalcGamma}, we have
\[y\in\widehat{\Ell}^k (\Z) \iff y\in\Z\ \text{and}\ \exists v\in \Lambda_k : y=p_2(v),\, p_1(v) \in W^k.\]
But if $y=p_2(v)$, then by the form of the matrix $M_{\lambda_1,\cdots,\lambda_k}$ we can write $v=(\widetilde v,0) + (0^{k-1},-y,y)$ with $\widetilde v\in \widetilde\Lambda_k$. Hence,
\begin{align*}
y\in\widehat{\Ell}^k (\Z) & \iff y\in\Z\ \text{and}\ \exists \widetilde v\in \widetilde\Lambda_k : (0^{k-1},-y)+\widetilde v	\in W^k\\
             & \iff y\in\Z\ \text{and}\ (0^{k-1},y)\in \bigcup_{\widetilde v\in\widetilde\Lambda_k} \widetilde v - W^k.
\end{align*}
Thus, $y\in\widehat{\Ell}^k (\Z)$ if and only if $y\in\Z$ and $\pr_{X_k}(0^{k-1},y)\in \pr_{X_k}(- W^k)$. Moreover, by construction, $\Varep_y^k = -w$.
\bigskip

Then, Proposition \ref{Flocker} follows directly from the fact that the points of the form $(0^{k-1},y)$, with $y\in\Z$, are equidistributed in $X_k$. To prove this equidistribution, we compute the inverse matrix of $\widetilde M_{\lambda_1,\cdots,\lambda_k}$:
\[{\widetilde M_{\lambda_1,\cdots,\lambda_k}}^{-1} = \begin{pmatrix}
\lambda_1^{-1} & \lambda_1^{-1}\lambda_2^{-1} &  \lambda_1^{-1}\lambda_2^{-1}\lambda_3^{-1} & \cdots & \lambda_1^{-1}\cdots \lambda_k^{-1}\\
    & \lambda_2^{-1}  & \lambda_2^{-1}\lambda_3^{-1} & \cdots & \lambda_2^{-1}\cdots \lambda_k^{-1}\\
    &      & \ddots & \ddots  & \vdots\\
    &      &        & \lambda_{k-1}^{-1} & \lambda_{k-1}^{-1}\lambda_k^{-1}\\
    &     &        &          & \lambda_k^{-1}
\end{pmatrix}.\]
Thus, the set of points of the form $(0^{k-1},y)$ in $X_k$ corresponds to the image of the map
\[\Z \ni y \longmapsto
{\widetilde M_{\lambda_1,\cdots,\lambda_k}}^{-1}\begin{pmatrix}
0^{k-1}\\ y \end{pmatrix}
=
\begin{pmatrix}
\widetilde\lambda_0^{-1}\\
\widetilde\lambda_1^{-1}\\
\vdots\\
\widetilde\lambda_{k-1}^{-1}\\
\end{pmatrix}y\]
in the canonical torus $\R^{k}/\Z^{k}$. But this map is ergodic when the family $( \widetilde\lambda_m^{-1})_{0\le m\le k}$ is $\Q$-free.
\end{proof}

From Proposition \ref{RoundoffLin2} it is possible to deduce an expression of the difference $c\delta(y)$ in terms of projections on a fundamental domain of $X_k$, as explained by the following proposition.

\begin{prop}\label{propDiscQuot}
The cumulated difference $c\delta(y)$ only depends on the projection of $(0^{k-1},y)$ on $X_k$. Moreover, the induced map $c\delta : X_k\to \R$ is affine when restricted to the fundamental domain
\[\D = \prod_{i=1}^k [1/2-\lambda_i,1/2]\]
of $X_k$, and if $(x_1,\dots,x_k) \in \D$ is the projection of $(0^{k-1},y)$ on $\D$ modulo $\widetilde\Lambda_k$, we have
\[c\delta(y) = - \frac12 - \sum_{m=1}^k x_m\frac{\widetilde\lambda_m}{\widetilde\lambda_0} .\]
\end{prop}

\begin{rem}\label{RemNormaliz}
From this proposition one can explain the appearance of normalisation constant in the definition of $c\delta$. Indeed, the mean of the affine map $c\delta$ on $\D$ is equal to its value in the centre of the parallelepiped $\D$:
\begin{align*}
\frac{1}{\widetilde\lambda_0} \int_{X_k} c\delta(x_1,\dots,x_k)\ud x_1\cdots\ud x_k & = c\delta\left(\frac12-\frac{\lambda_1}{2}\,,\ \cdots\,,\ \frac12-\frac{\lambda_k}{2}\right)\\
      & = -\frac12 - \sum_{m=1}^k\widetilde\lambda_m\frac{1-\lambda_m}{2\widetilde\lambda_0}\\
	  & = -\frac{1}{2\widetilde\lambda_0}.
\end{align*}
In particular, for $R\in\N$, one has
\begin{align*}
\frac{1}{R}\int_0^R c\delta(y)\ud y & = \frac{1}{R} \sum_{n=0}^{R-1}c\delta(n+1/2)\\
  & = \frac{1}{R} \sum_{n=0}^{R-1} \left(c\delta(n) + \frac{1}{2\widetilde\lambda_0}\right)\\
  & \underset{R\to +\infty}{\longrightarrow} 0,
\end{align*}
in other words the map $c\delta$ has zero mean.
\end{rem}

\begin{proof}[Proof of Proposition \ref{propDiscQuot}]
Given $n\in\N$, we want to compute the cumulated difference \eqref{DefDiscR}:
\[c\delta(n) = \frac{n}{\widetilde\lambda_0} - \card\big\{x\in\N\mid \widehat{\Ell}^k (x)\le n\big\} + \frac{1}{2}\]
In other words (as $x\mapsto \widehat{\Ell}^k(x)$ is increasing) we search for the biggest $x\in\N$ such that $\widehat{\Ell}^k(x) \le n$. Let $x$ be such a number, we have
\begin{equation}\label{EqPropCDelta}
c\delta(n) = \frac{n}{\widetilde\lambda_0} - x - \frac{1}{2}
\end{equation}
(remark that this formula allows us to read this Cramér distance on the ``time 0'' set $\Z$ -- see Figure \ref{FigDisc} --, this is possible by the preservation of order of the maps $\ell_i$ and $\widehat \ell_i$). We are reduced to compute this integer $x$.

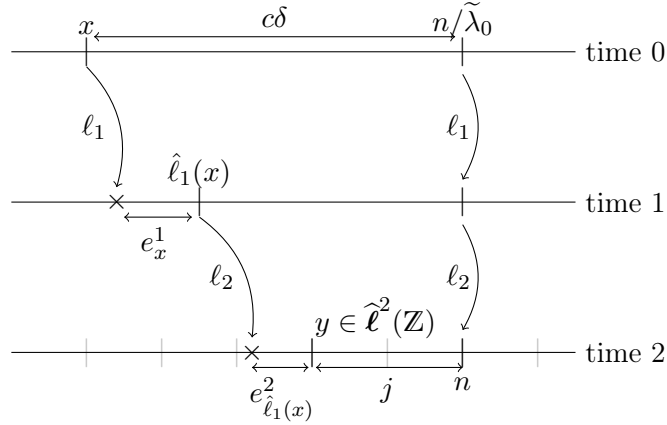
\begin{figure}
\begin{center}
\begin{tikzpicture}[scale=1]

\draw(1,0) -- (8.5,0);
\draw (8.5,0) node[right]{time 0};
\draw(1,-2) -- (8.5,-2);
\draw (8.5,-2) node[right]{time 1};
\foreach\i in {2,...,8}{\draw[color=gray!50!white] (\i,-4) node{$|$};}
\draw(1,-4) -- (8.5,-4);
\draw (8.5,-4) node[right]{time 2};
\draw(2,0) node{$|$};
\draw(2.4,-2) node{$\times$};
\draw(3.5,-2) node{$|$};
\draw(4.2,-4) node{$\times$};
\draw(5,-4) node{$|$};
\draw(7,-4) node{$|$};
\draw(7,-2) node{$|$};
\draw(7,0) node{$|$};

\draw[->](2,-.2) to[bend left] (2.4,-1.8);
\draw[<->](2.5,-2.2) -- (3.4,-2.2);
\draw[->](3.5,-2.2) to [bend left] (4.2,-3.8);
\draw[<->](4.2,-4.2) -- (4.95,-4.2);
\draw[<->](5.05,-4.2) -- (7,-4.2);
\draw[<->](2.1,.2) -- (6.9,.2);
\draw[->](7,-2.3) to [bend left] (7,-3.7);
\draw[->](7,-.3) to [bend left] (7,-1.7);

\draw(2,0.1) node[above] {$x$};
\draw(2.9,-2.2) node[below] {$e_x^1$};
\draw(3.5,-1.95) node[above] {$\hat \ell_1(x)$};
\draw(2.1,-1) node {$\ell_1$};
\draw(6.95,-1) node {$\ell_1$};
\draw(3.8,-3) node {$\ell_2$};
\draw(6.95,-3) node {$\ell_2$};
\draw(4.6,-4.2) node[below] {$e_{\hat\ell_1(x)}^2$};
\draw(6,-4.2) node[below] {$j$};
\draw(4.9,-3.95) node[above right] {$y\in\widehat\Ell^2(\Z)$};
\draw(7,-4.2) node[below] {$n$};
\draw(7,0.05) node[above] {$n/\widetilde\lambda_0$};
\draw(4.5,0.2) node[above] {$c\delta$};

\end{tikzpicture}
\end{center}
\caption{Computation of the cumulated difference $c\delta$ in the case $k=2$. The point $y$ is the biggest point of $\Ell^2(\Z)$ smaller than $n$.}\label{FigDisc}
\end{figure}

We denote $y = \widehat{\Ell}^k(x)$ and $j = n-y\in\N$. 
In this case, the definition of global error leads to
\[\mathcal E_x^k = \widehat\Ell^k(x) - \Ell^k(x) = y - \widetilde\lambda_0 x \qquad \iff \qquad  x = \frac{1}{\widetilde \lambda_0}\big( y-\mathcal E_x^k\big).\]
Thus, applying this to \eqref{EqPropCDelta},
\[c\delta(n) = \frac{n}{\widetilde\lambda_0} - \frac{y-\mathcal E_x^k}{\widetilde\lambda_0} - \frac12 = \frac{j + \mathcal E_x^k}{\widetilde\lambda_0} - \frac12.\]
But, combining Formula \eqref{FormulEx} page \pageref{FormulEx} linking the global error $\mathcal E_x$ with the roundoff error vector $\Varep_x^k$ with the fact that $\Varep_x^k = -w$ (Proposition \ref{RoundoffLin2}), one has
\[\mathcal E_x^k = -\sum_{m=1}^k \widetilde\lambda_m w_m,\]
so that
\[c\delta(n) = \frac{1}{\widetilde\lambda_0} \Big(j - \sum_{m=1}^k \widetilde\lambda_m w_m\Big) - \frac12,\]
with $w$ depending only on the projection of $(0^{k-1},y)=(0^{k-1},n-j)$ on $X_k$.

We have reduced to find, given $n\in\Z$, the smallest $j\in\N$ such that $n-j \in \widehat{\Ell}^k(\Z)$. First remark that $W^k$ projects injectively (but not surjectively) on the torus $X_k = \R^k/\widetilde \Lambda_k$. So we define a partition of $X_k$ into first visit sets $W_j$ in $W^k$ under the iterates of the translation $(0^{k-1},1)$ (see Figure \ref{figPartition}). More precisely, $v\in X_k$ belongs to $W_j$ iff $j$ is the smallest nonegative integer such that $v + (0^{k-1},j) \in W^k$ mod $\widetilde \Lambda_k$.

Thus, if $(0^{k-1},-n)\in W_j$, then $j$ is the smallest integer such that $(0^{k-1},-(n-j))\in W_0 = W^k$. In this case, $j$ is the smallest nonnegative integer such that $n-j\in \widehat\Ell^k(\Z)$, moreover
\[W_i = \pr_{X_k}\big(W^k - (0^{k-1},i)\big) \setminus \bigcup_{j=0}^{i-1} W_j.\]

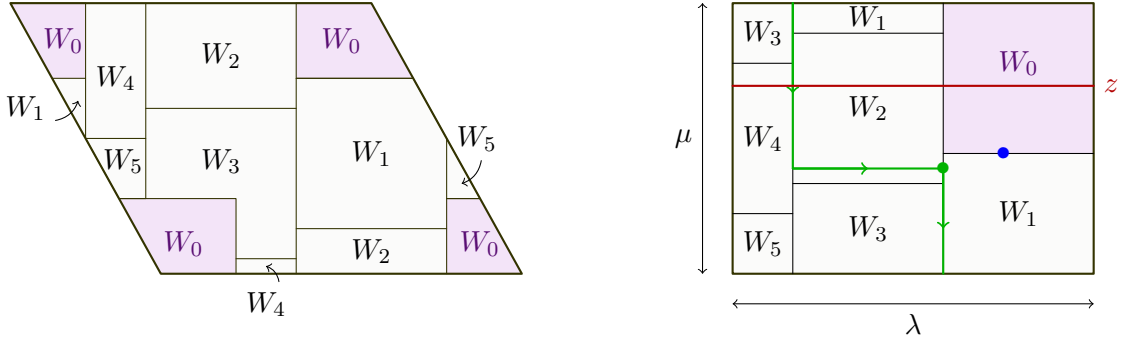
\begin{figure}
\begin{center}
\definecolor{zzttqq}{rgb}{0.2,0.2,0.}
\definecolor{caca}{rgb}{0.7,0.1,0.9}
\begin{tikzpicture}[line cap=round,line join=round, scale=2, baseline]

\fill[color=zzttqq,fill opacity=0.02] (0.,0.) -- (2.4,0.) -- (1.4,1.8) -- (-1.,1.8) -- cycle;
\fill[color=zzttqq,fill=caca,fill opacity=0.1] (-0.2777777777777778,0.5) -- (0.5,0.5) -- (0.5,0.) -- (0.,0.) -- cycle;
\fill[color=zzttqq,fill=caca,fill opacity=0.1] (1.9,0.5) -- (2.1222222222222222,0.5) -- (2.4,0.) -- (1.9,0.) -- cycle;
\fill[color=zzttqq,fill=caca,fill opacity=0.1] (1.4,1.8) -- (1.677777777777778,1.3) -- (0.9,1.3) -- (0.9,1.8) -- cycle;
\fill[color=zzttqq,fill=caca,fill opacity=0.1] (-1.,1.8) -- (-0.5,1.8) -- (-0.5,1.3) -- (-0.7222222222222222,1.3) -- cycle;

\draw [color=zzttqq, thick] (0.,0.)-- (2.4,0.)-- (1.4,1.8) -- (-1.,1.8) -- cycle;
\draw [color=zzttqq] (-0.2777777777777778,0.5) -- (0.5,0.5)--(0.5,0.);
\draw [color=zzttqq] (1.9,0.5)-- (2.1222222222222222,0.5);
\draw [color=zzttqq] (1.9,0.)-- (1.9,0.9);
\draw [color=zzttqq] (1.677777777777778,1.3)-- (0.9,1.3);
\draw [color=zzttqq] (0.9,0)-- (0.9,1.8);
\draw [color=zzttqq] (0.9,1.8)-- (1.4,1.8);
\draw [color=zzttqq] (-1.,1.8)-- (-0.5,1.8);
\draw [color=zzttqq] (-0.5,1.8)-- (-0.5,0.9);
\draw [color=zzttqq] (-0.5,1.3)-- (-0.7222222222222222,1.3);
\draw [color=zzttqq] (0.9,1.3)-- (1.677777777777778,1.3);
\draw [color=zzttqq] (1.677777777777778,1.3)-- (1.9,0.9);
\draw [color=zzttqq] (1.9,0.3)-- (0.9,0.3);
\draw [color=zzttqq] (-0.1,1.8)-- (-0.1,0.5);
\draw [color=zzttqq] (-0.1,1.1)-- (0.9,1.1);
\draw [color=zzttqq] (0.5,0.1)-- (0.9,0.1);
\draw [color=zzttqq] (-0.1,0.9)-- (-0.5,0.9);
\draw [color=zzttqq] (-0.1,0.5)-- (-0.2777777777777778,0.5);
\draw [color=zzttqq] (-0.2777777777777778,0.5)-- (-0.5,0.9);
\draw [color=zzttqq] (1.9,0.5)-- (2.1222222222222222,0.5);
\draw [color=zzttqq] (2.1222222222222222,0.5)-- (1.9,0.9);

\node[color=caca!50!black] at (1.2,1.55) {$W_0$};
\node[color=caca!50!black] at (0.15,0.2) {$W_0$};
\node[color=caca!50!black] at (2.12,0.2) {$W_0$};
\node[color=caca!50!black] at (-.65,1.55) {$W_0$};
\node  at(1.4,0.8) {$W_1$};
\node  at(1.4,0.15) {$W_2$};
\node  at(0.4,1.45) {$W_2$};
\node  at(0.4,0.75) {$W_3$};
\node  at(-.3,1.3) {$W_4$};
\node  at(-.25,0.75) {$W_5$};
\node(A) at(-.9,1.1) {$W_1$};
\draw[->] (A) to[bend right] (-.55,1.15);
\node(B) at(2.1,0.9) {$W_5$};
\draw[->] (B) to[bend left] (2,0.6);
\node(C) at(0.7,-.2) {$W_4$};
\draw[->] (C) to[bend right] (0.7,0.05);

\end{tikzpicture}\hfill
\begin{tikzpicture}[line cap=round,line join=round, scale=2, baseline]
\fill[color=zzttqq,fill opacity=0.02] (0.,0.) -- (2.4,0.) -- (2.4,1.8) -- (0.,1.8) -- cycle;
\fill[color=zzttqq,fill=caca,fill opacity=0.1] (2.4,1.8) -- (1.4,1.8) -- (1.4,0.8) -- (2.4,0.8)  -- cycle;
\draw [color=zzttqq, thick] (0.,0.) -- (2.4,0.) -- (2.4,1.8) -- (0.,1.8) -- cycle;
\draw(1.4,1.8) -- (1.4,0);
\draw(0.4,1.8) -- (0.4,0);
\draw(1.4,0.8) -- (2.4,0.8);
\draw(0.4,1.6) -- (1.4,1.6);
\draw(0.4,0.6) -- (1.4,0.6);
\draw(0,1.4) -- (0.4,1.4);
\draw(0,0.4) -- (0.4,0.4);
\draw(1.8,0.8) node[color=blue]{$\bullet$};
\draw(1.4,0.7) node[color=green!70!black]{$\bullet$};
\draw[color=green!70!black, thick](0.4,1.8) -- (0.4,0.7) -- (1.4,0.7) -- (1.4,0);
\draw[color=green!70!black, thick, ->](1.4,0.7) -- (1.4,0.3);
\draw[color=green!70!black, thick, ->](0.4,1.8) -- (0.4,1.2);
\draw[color=green!70!black, thick, ->](0.4,0.7) -- (0.9,0.7);

\node[color=caca!50!black] at (1.9,1.4) {$W_0$};
\node  at(1.9,0.4) {$W_1$};
\node  at(0.9,1.7) {$W_1$};
\node  at(0.9,1.1) {$W_2$};
\node  at(0.9,0.3) {$W_3$};
\node  at(0.2,1.6) {$W_3$};
\node  at(0.2,0.9) {$W_4$};
\node  at(0.2,0.2) {$W_5$};

\draw[<->] (0,-.2) -- (2.4,-.2) node[midway, below]{$\lambda$};
\draw[<->] (-.2,1.8) -- (-.2,0) node[midway, left]{$\mu$};

\draw[color=red!70!black, thick](0,1.25) -- (2.4,1.25);
\draw[color=red!70!black] (2.4,1.25) node[right] {$z$};

\end{tikzpicture}
\end{center}
\vspace{-15pt}
\caption{\label{figPartition}The sets $W_i$ for $\lambda=2.4$ and $\mu=1.8$ on a canonical fundamental domain (left) and a rectangular fundamental domain $\mathcal D$ (right). We recognize a suspension of the rotation $x\mapsto x-1$ modulo $\lambda$.}
\end{figure}

For $v\in X_k$, this allows to define
\[c\delta(v) = \frac{1}{\widetilde\lambda_0} \Big(j - \sum_{m=1}^k \widetilde\lambda_m w_m\Big) - \frac12,\]
where $w = v + (0^{k-1},j) \in W^k$ (and hence $v\in W_j$).

We just have proved that the map $c\delta$ is piecewise affine on $X_k$ (more precisely, affine in restriction to each set $W_i$). Let us analyse the partition into the sets $W_j$ a bit further by looking at its projection on the fundamental domain $\D = \prod_{i=1}^k [1/2,1/2-\lambda_i]$ (see Figure \ref{figPartition}). In $X_k$, the set $W_{j+1}$ is simply obtained by a translation of $W_j$ by $(0,-1)$, removing the intersections with the $W_i$ for $i\le j$ if necessary. In $\D$, the gluings of opposite faces made to recover $X_k$ correspond to suspensions of rotations $x\mapsto x+e_i$ where $e_i$ is the $i$-th vector of the canonical basis of $\R^k$. In particular, one sees that all the $W_i$, except from the  $\lceil \lambda_k\rceil$ last, are simply translates of $W_0$. 

It remains to prove that the map $v\mapsto c\delta(v)$ is continuous on $\D$. First remark that it is continuous (because linear) on every set $W_j$, so we just have to prove the continuity at the boundaries of the $W_i$'s. To do it we reason by recurrence on the dimension.

As a first step, we prove that $c\delta$ is continuous in restriction to the union of two cubes having a face in common orthogonal to the last canonical coordinate vector. In this case, they are two cubes with consecutive indices, say $j$ and $j+1$. Let us take a point $v$ belonging to their common face (the blue point of Figure \ref{figPartition}, right). On the one hand, $v\in W_j$, so we can write $v = (w_1,\cdots,w_k)$, where the $w_i$'s are the coordinates of $v$ with respect to the centre of the set $W_i$. This leads to
\[c\delta(v) = \frac{1}{\widetilde\lambda_0} \Big(j - \sum_{m=1}^{k-1} \widetilde\lambda_m w_m - \widetilde\lambda_k w_k\Big) - \frac12.\]
But $w_k = -1/2$ and $\widetilde\lambda_k=1$, so 
\[c\delta(v) = \frac{1}{\widetilde\lambda_0} \Big(j - \sum_{m=1}^{k-1} \widetilde\lambda_m w_m +\frac12 \Big) - \frac12.\]
On the other hand, we also have $v\in W_{j+1}$, and so if we write $v = (w'_1,\cdots,w'_k)$ in the coordinates with respect to the centre of $W_{j+1}$,then 
\[c\delta(v) = \frac{1}{\widetilde\lambda_0} \Big(j+1 - \sum_{m=1}^{k-1} \widetilde\lambda_m w'_m - \widetilde\lambda_k w'_k\Big) - \frac12.\]
But we have $w_i = w'_i$ for all $i\le k-1$ and $w'_k = 1/2$, so
\[c\delta(v) = \frac{1}{\widetilde\lambda_0} \Big(j+1 - \sum_{m=1}^{k-1} \widetilde\lambda_m w_m -\frac12 \Big) - \frac12.\]
We deduce that both values of $c\delta(v)$ coincide whether $v$ is seen a an element of $W_j$ or of $W_{j+1}$.

For the induction's heredity, we use the trivial formula
\[\sum_{m=1}^{k-1}\widetilde\lambda_m w_m = \sum_{m=1}^{k-1}\widetilde\lambda_m w'_m,\]
where $w_m - w'_m$ is a vector having zero coordinates but the $\ell$-th one equal to $-1$ and the $\ell+1$-th one equal to $\lambda_{\ell+1}$. For each $\ell$, using the fact that $c\delta$ is continuous with respect to the $\ell+1$-th coordinate, this tells us that the map $c\delta$ is continuous with respect to the $\ell$-th coordinate (geometrically, it consists in following the green path of Figure \ref{figPartition}, right).

Remark that one can also prove the continuity by examinating directly what happens on the image spaces $\Ell^i(\R)$ and performing small translations of the grids.

\end{proof}

As a corollary, one can get an alternative proof of Proposition \ref{WasDistLin}.

\begin{proof}[Second proof of Proposition \ref{WasDistLin}]\label{Proof2}
By Equation \eqref{EqCalcDiscdisc} page \pageref{EqCalcDiscdisc}, one has
\[\Disc^2 = \frac{1}{12 \widetilde\lambda_0^2} + \Var\left( c\delta\Big(n+\frac12\Big)\right).\]
The second term corresponds to the variance of the map $c\delta + 1/(2\widetilde\lambda_0)$ on $\D$. As this map is affine with zero mean, and by the form of $\D$, this variance is equal to the sum of the variances of its coordinates, i.e.
\begin{align*}
\Var \left(c\delta\Big(n+\frac12\Big)\right) & = \sum_{m=1}^k \Var_{[1/2-\lambda_m,1/2]}\left(x_m \frac{\widetilde\lambda_m}{\widetilde\lambda_0}\right)\\
    & = \sum_{m=1}^k \frac{\lambda_m^2\widetilde\lambda_m^2}{12\widetilde\lambda_0^2},
\end{align*}
thus
\[\Disc = \frac{1}{12\widetilde\lambda_0^2}\sum_{m=0}^k \widetilde\lambda_{m}^2.\]
\end{proof}

\section{The tree linear case}\label{Sectree}

In this section we adapt the study made in the previous one to stick to the shape of the set of preimages of a point under some expanding map, which has a structure of $d$-ary tree. We will get get Cramér distance estimations for a complete $d$-ary tree with edges decorated by linear expanding maps. 

Fix $r\ge 1$ and $d\ge 2$. We begin by the definition of the set of expanding maps.

\begin{definition}\label{DefExpan}
We denote by $\mathcal D^r(\Sp^1)$ the set of $C^r$ expanding maps of degree $d$ of $\Sp^1$. More precisely, $\mathcal D^r(\Sp^1)$ is the set of degree $d$ maps $f : \Sp^1\to \Sp^1$ such that the derivative $f^{(\lfloor r\rfloor)}$ is well defined and belongs to $C^{r-\lfloor r\rfloor}(\T^n)$ and such that for every $x\in \Sp^1$, we have $|f'(x)| > 1$.
\end{definition}

The set of preimages of a point $x\in \Sp^1$ by an expanding map $f$ has a natural structure of complete $d$-ary tree. We now define the linear setting corresponding to the local behaviour of $f\in\mathcal D^r(\Sp^1)$ using this representation. 

\begin{definition}
We set (see also Figure~\ref{TreeBee})\index{$I_k$}
\[I_k = \bigsqcup_{m=1}^k \llbracket 1,d\rrbracket^m\]
the set of $m$-tuples of integers of $\llbracket 1,d\rrbracket$, for $1\le m\le k$.

For $\ind = (i_1,\cdots,i_m)\in \llbracket 1,d\rrbracket^m$, we set $\len(\ind) = m$\index{$\len(\ind)$} and the parent $\fat(\ind) = (i_1,\cdots,i_{m-1})\in \llbracket 1,d\rrbracket^{m-1}$\index{$\fat(\ind)$} (with the convention $\fat(i_1) = \emptyset$).
\end{definition}

\begin{figure}[!b]
\begin{minipage}[c]{.39\linewidth}
\begin{center}
\begin{tikzpicture}[scale=.7]
\node (O) at (1,0){$\emptyset$};
\node (A) at (3,1){$(1)$};
\node (B) at (3,-1){$(2)$};
\node (C) at (6,1.5){$(1,1)$};
\node (D) at (6,.5){$(1,2)$};
\node (E) at (6,-.5){$(2,1)$};
\node (F) at (6,-1.5){$(2,2)$};
\draw (O) -- (A);
\draw (A) -- (C);
\draw (A) -- (D);
\draw (O) -- (B);
\draw (B) -- (E);
\draw (B) -- (F);
\end{tikzpicture}
\caption[The tree $T_2$ for $d=2$]{The tree $T_2$ for $d=2$.}\label{TreeBee}
\end{center}
\end{minipage}\hfill
\begin{minipage}[c]{.6\linewidth}
\begin{center}
\begin{tikzpicture}[scale=.7]
\node (O) at (1,0){$y$};
\node (A) at (3,1){$x_{(1)}$};
\node (B) at (3,-1){$x_{(2)}$};
\node (C) at (6,1.5){$x_{(1,1)}$};
\node (D) at (6,.5){$x_{(1,2)}$};
\node (E) at (6,-.5){$x_{(2,1)}$};
\node (F) at (6,-1.5){$x_{(2,2)}$};
\draw (O) -- (A);
\draw (A) -- (C);
\draw (A) -- (D);
\draw (O) -- (B);
\draw (B) -- (E);
\draw (B) -- (F);
\end{tikzpicture}
\caption[Tree associated to the preimages of $y$]{The tree associated to the preimages of $y$, for $k=2$ and $d=2$. We have $f(x_{(1,1)}) = f(x_{(1,2)}) = x_{(1)}$, etc.}\label{ProbTree}
\end{center}
\end{minipage}
\end{figure}

The set $I_k$ is the linear counterpart of the set $\bigsqcup_{m=1}^k f^{-m}(\{y\})$. Its cardinal is equal to $d(1-d^k)/(1-d)$.

\begin{definition}\label{Noel!}
Let $k\in\N$. The \emph{complete tree of order $k$} is the rooted $d$-ary tree $T_k$\index{$T_k$} whose vertices are the elements of $I_k$ together with the root $\emptyset$, and whose edges are of the form $(\fat(\ind),\ind)_{\ind\in I_k}$ (see Figure~\ref{TreeBee}).
%
%
%
%
\end{definition}

We now consider  a family $(\ell_\ind)_{\ind\in I_k}$ of homotheties of $\R$ of parameters $(\lambda_\ind)_{\ind\in I_k}>1$. In this new case we denote
\[\widehat{\Ell}^k\big(\Z^{d^k}\big) = \bigcup_{\ind\in \llbracket 1,d\rrbracket^k}\big( \widehat \ell_{\fat^{k-1}(\ind)} \circ \cdots \circ \widehat \ell_\ind \big) (\Z),\]
\[\widetilde \lambda_{\ind} = \lambda_{\ind}\lambda_{\fat(\ind)} \cdots \lambda_{\fat^{\len(\ind)-1}(\ind)}, \qquad  \widetilde\lambda_{tot}^{-1} = \sum_{\ind\in \llbracket 1,d\rrbracket^k}\widetilde\lambda_\ind^{-1},\]
and for $\ind\in\llbracket 1,d\rrbracket^k$ (see also \eqref{DefDiscR}, we omit the dependance on the measures),
\begin{equation}\label{DefDisc2}
c\delta_\ind(y) = \frac{y}{\widetilde \lambda_\ind} - \card\big\{x\in\N\mid \widehat \ell_{\fat^{k-1}(\ind)} \circ \cdots \circ \widehat \ell_\ind (x)\le y\big\} + \frac{1}{2}
\end{equation}
and, denoting $\Er_{x,\ind}^k$ the cumulative error made by iterating the point $x$ following the path from $\ind$ during time $k$. 
\begin{align*}
\Disc_{R}^2
= & \frac{1}{R} \int_0^{R} \left(\widetilde\lambda_{tot}^{-1} y - \sum_{\ind\in \llbracket 1,d\rrbracket^k}\Big(\max\big\{x\mid \widetilde\lambda_\ind x+\Er_{x,\ind}^k \le y\big\} - \frac12\Big) \right)^2\ud y \\
= & \frac{1}{R} \int_0^{R} \left(\sum_{\ind\in \llbracket 1,d\rrbracket^k} c\delta_\ind(y) \right)^2\ud y.
\end{align*}

For $\ind, \ind'\in I_k$ of same length, we denote
\[k_0(\ind,\ind') = \min\big\{m\in \llbracket 0,k\rrbracket \mid \fat^m(\ind) = \fat^m(\ind')\big\}.\]
We also denote $\delta_{\Z^{d^k}} = \sum_{x\in\Z^{d^k}}\delta_x$ the uniform measure on $\Z^{d^k}$

\begin{prop}\label{WasDistTree}
Let $k\in\N$, and a family $(\ell_\ind)_{\ind\in I_k}$ of homotheties of $\R$ of parameters $(\lambda_\ind)_{\ind\in I_k}$ strictly bigger than 1. If the family $(\widetilde\lambda_{\ind})_{\ind\in I_k }$ is $\Q$-free (which is a generic condition), then
\[
\Disc^2 \left(\widehat{\Ell}^k_*\big(\delta_{\Z^{d^k}}\big)\, , \, \widetilde\lambda_{tot}^{-1} \Leb\right)
  = \frac{1}{12\widetilde\lambda_{tot}} + \frac{1}{12} \sum_{\ind,\ind'\in \llbracket 1,d\rrbracket^k} \sum_{m=k_0(\ind,\ind')}^{k-1} \frac{\widetilde\lambda_{\fat^m(\ind)}\widetilde\lambda_{\fat^m(\ind')}}{\widetilde\lambda_{\ind}\widetilde\lambda_{\ind'}}.
\]

\end{prop}

Note that for any $m \ge k_0(\ind,\ind')$, one has $\widetilde\lambda_{\fat^m(\ind)} = \widetilde\lambda_{\fat^m(\ind')}$.

To prove this proposition we will need to estimate the expectation of the following correlations:
\[\Corr_{\ind,\ind'}(y) = c\delta_\ind\left(y+\frac12\right) c\delta_{\ind'}\left(y+\frac12\right).\]

\begin{lemme}\label{lemEstimCorr}
Under the hypotheses of Proposition~\ref{WasDistTree}, for any $\ind\neq \ind'$, $\mathbb E\big[\Corr_{\ind,\ind'}\big]$ is well defined (that is, the limit exists) and satisfies
\[\mathbb E\big[\Corr_{\ind,\ind'}(y)\big] = \frac{1}{12 \widetilde\lambda_\ind\widetilde\lambda_{\ind'}} \sum_{m=k_0(\ind,\ind')}^{k-1} \widetilde\lambda_{\fat^m(\ind)}\widetilde\lambda_{\fat^m(\ind')} \ge 0.\]
\end{lemme}

This lemma is the heart of the proof of Proposition \ref{WasDistTree}. It is deduced from the study conducted in the previous section (Proposition \ref{propDiscQuot}). We first deduce Proposition \ref{WasDistTree} from it.

\begin{proof}[Proof of Proposition \ref{WasDistTree}]
By Lemma \ref{LemCalcDiscdisc}, for any $R\in\N$,
\[\Disc_{R}^2 = \frac{1}{R} \sum_{y=0}^{R-1} \left(\sum_{\ind\in \llbracket 1,d\rrbracket^k} c\delta_\ind\left(y+\frac12\right) \right)^2 + \frac{1}{12\widetilde\lambda_{tot}}.\]
But
\[
\left(\sum_{\ind\in \llbracket 1,d\rrbracket^k} c\delta_\ind\left(y+\frac12\right) \right)^2 = \sum_{\ind\in \llbracket 1,d\rrbracket^k} c\delta_\ind\left(y+\frac12\right)^2 + \sum_{\ind\neq\ind'\in \llbracket 1,d\rrbracket^k} c\delta_\ind\left(y+\frac12\right) c\delta_{\ind'}\left(y+\frac12\right).
\]
By Corollary~\ref{CoroWasDistLin}, the first term has mean
\[\frac{1}{12} \sum_{\ind\in \llbracket 1,d\rrbracket^k} \frac{1}{\widetilde\lambda_\ind^2}\sum_{m=0}^{k-1} \widetilde\lambda_{\fat^m(\ind)}^2,\]
while that of the second one is given by Lemma \ref{lemEstimCorr}.
\end{proof}

\begin{proof}[Proof of Lemma \ref{lemEstimCorr}]
Let $\ind,\ind' \in\llbracket 1,d\rrbracket^k$ such that $\ind\neq\ind'$. We write $\ind = (i_k,\dots,i_1)$ and $\ind' = (i'_k,\dots,i'_1)$ (note the fact that the indices decrease, to correspond to the order of iterations of maps as in the previous part). We also note $\widetilde\Lambda_\ind$ the lattice associated to the multipliers $\lambda_{i_1},\dots,\lambda_{i_k}$, see \eqref{DefMat} page \pageref{DefMat} (and the same for $\widetilde\Lambda_{\ind'}$).

By Proposition \ref{propDiscQuot}, we know that for $y\in\Z$, the cumulated difference $c\delta_\ind(y)$ only depends on the projection of $(0^{k-1},y)$ on the torus $\R^k/\widetilde \Lambda_\ind$ and is affine when restricted to the fundamental domain $\D = \prod_{j=1}^k [1/2,1/2-\lambda_{i_j}]$ of $\widetilde\Lambda_\ind$: if $(0^{k-1},y) = (x_1,\dots,x_k) \in \D \mod \widetilde\Lambda_\ind$, we have
\[c\delta_\ind(y+1/2) = - \frac12 - \sum_{m=1}^k x_m\frac{\widetilde\lambda_{\fat^{m-1}(\ind)}}{\widetilde\lambda_\ind} + \frac{1}{2\widetilde\lambda_{\ind}}.\]
Thus, for any $y\in\Z$ (with transparent notations),
\[\Corr_{\ind,\ind'}(y) = \widetilde\lambda_\ind^{-1}\widetilde\lambda_{\ind'}^{-1}
\left(\frac{\widetilde\lambda_\ind-1}{2} + \sum_{m=1}^k x_m \widetilde\lambda_{\fat^{m-1}(\ind)}\right)
\left(\frac{\widetilde\lambda_{\ind'}-1}{2} + \sum_{m=1}^k x'_m \widetilde\lambda_{\fat^{m-1}(\ind')}\right).\]
By definition of $k_0\overset{\text{def.}}{=} k_0(\ind,\ind')$, for every $m\in\llbracket k_0+1,k\rrbracket$, we have $\widetilde\lambda_{\fat^{m-1}(\ind)} = \widetilde\lambda_{\fat^{m-1}(\ind')}$ and $x_m = x'_m$; thus we can denote
\[t_1 = \sum_{m=k_0+1}^k x_m \widetilde\lambda_{\fat^{m-1}(\ind)} =
\sum_{m=k_0+1}^k x'_m\widetilde\lambda_{\fat^{m-1}(\ind')},\]
\[t_2 = \frac{\widetilde\lambda_{\ind}-1}{2} + \sum_{m=1}^{k_0} x_m \widetilde\lambda_{\fat^{m-1}(\ind)} \qquad \text{and} \qquad t'_2 = \frac{\widetilde\lambda_{\ind'}-1}{2} + \sum_{m=1}^{k_0} x'_m \widetilde\lambda_{\fat^{m-1}(\ind')},\]
so that
\[\widetilde\lambda_\ind\widetilde\lambda_{\ind'} \Corr_{\ind,\ind'}(y) = \big(t_1 + t_2\big)\big(t_1+t'_2\big).\]
Thus
\begin{align*}
\widetilde\lambda_\ind\widetilde\lambda_{\ind'} \mathbb E\big[\Corr_{\ind,\ind'}(y)\big] = &
\mathbb E[t_1^2] + \mathbb E[t_1t_2'] + \mathbb E[t_1t_2] + \mathbb E[t_2t_2']\\
 = & \mathbb E[(t_1-\mathbb E[t_1])^2] + \mathbb E[t_1]^2 + \mathbb E[t_1t_2']
+ \mathbb E[t_1t_2] + \mathbb E[t_2t_2'].
\end{align*}
The hypothesis of independence over $\Q$ of the family $(\widetilde\lambda_{\ind})_{\ind\in I_k }$ implies that the events $t_1$, $t_2$ and $t_2'$ are ``independent''. To see it, consider the (big) matrix \setlength{\arraycolsep}{3pt}
\begin{multline*}
M_{k+k_0}(\R) \ni \widetilde M_{\ind,\ind'} = \\
\left(\begin{array}{cccc|cccc|ccc}
\lambda_\ind & -1 &   &   &   &   &   &   & \\
    & \ddots & \ddots &   &   &   &   &   &\\
    &   & \lambda_{\fat^{k_0-2}(\ind)} & -1 &   &   &   & \\
    &   &   & \lambda_{\fat^{k_0-1}(\ind)} & 0 &   &   &   & -1\\ \hline 
    &   &   &   & \lambda_{\ind'}  & -1 &   & \\
    &   &   &   &   & \ddots & \ddots &   & \\
    &   &   &   &   &   & \lambda_{\fat^{k_0-2}(\ind')} & -1 \\
    &   &   &   &   &   &   & \lambda_{\fat^{k_0-1}(\ind')} & -1 \\ \hline
    &   &   &   &   &   &   &   & \lambda_{\fat^{k_0}(\ind)} & \ddots \\
    &   &   &   &   &   &   &   &   & \ddots & -1\\
    &   &   &   &   &   &   &   &   &   & \lambda_{\fat^{k-1}(\ind)}
\end{array}\right),
\end{multline*}
remark that the quotient space $\R^{k+k_0}/ \widetilde M_{\ind,\ind'}\Z^{k+k_0}$ is obtained from the product space $\R^k / \widetilde M_{\lambda_{i_1},\cdots,\lambda_{i_k}}\Z^k \times \R^k / \widetilde M_{\lambda_{i'_1},\cdots,\lambda_{i'_k}}\Z^k$ by quotienting by the $k-k_0$ last coordinates, which corresponds to identical multipliers $\lambda_{\fat^j(\ind)} = \lambda_{\fat^j(\ind')}$, thus to identical projections of $(0^{k-1},y)$ on the fundamental domains $\D$ and $\D'$. A simple computation of $\widetilde M_{\ind,\ind'}^{-1}$, similar to that of the proof of Proposition \ref{RoundoffLin2}, leads to the fact that when the family $(\widetilde\lambda_{\ind})_{\ind\in I_k }$ is $\Q$-free, then
\begin{align*}
y\mapsto \widetilde M_{\ind,\ind'}^{-1}(0,y) = y \Big( & \widetilde\lambda_\ind^{-1},\ 
\widetilde\lambda_{\fat(\ind)}^{-1},\ 
\cdots,\ 
\widetilde\lambda_{\fat^{k_0-1}(\ind)}^{-1},\\
& \widetilde\lambda_{\ind'}^{-1},\ 
\widetilde\lambda_{\fat(\ind')}^{-1},\ 
\cdots,\ 
\widetilde\lambda_{\fat^{k_0-1}(\ind')}^{-1},\\
& \widetilde\lambda_{\fat^{k_0}(\ind)}^{-1},\ 
\widetilde\lambda_{\fat^{k_0+1}(\ind)}^{-1},\ 
\cdots,\ 
\widetilde\lambda_{\fat^{k-1}(\ind)}^{-1}\Big)
\end{align*}
is an ergodic $\Z$-action on the torus $\R^{k+k_0}/\Z^{k+k_0}$. In particular,
\[\mathbb{E}[t_1t_2] = \mathbb{E}[t_1]\mathbb{E}[t_2], \quad \mathbb{E}[t_1t_2'] = \mathbb{E}[t_1]\mathbb{E}[t_2'] \quad \text{and}\quad \mathbb{E}[t_2t_2'] = \mathbb{E}[t_2]\mathbb{E}[t_2'].\]

\bigskip

We deduce that
\begin{align*}
\widetilde\lambda_\ind\widetilde\lambda_{\ind'}\mathbb E\big[\Corr_{\ind,\ind'}(y)\big] = & \mathbb E\big[(t_1-\mathbb E[t_1])^2\big] + \mathbb E[t_1]^2 + \mathbb E[t_1]\mathbb E[t_2'] + \mathbb E[t_1]\mathbb E[t_2] + \mathbb E[t_2]\mathbb E[t_2']\\
= & \mathbb E\big[(t_1-\mathbb E[t_1])^2\big] + \mathbb E\big[t_1 + t_2\big] \mathbb E\big[t_1+t'_2\big].
\end{align*}
But by Remark \ref{RemNormaliz}, $\mathbb E[t_1 + t_2] \mathbb E[t_1+t'_2] = 0$, so
\[\widetilde\lambda_\ind\widetilde\lambda_{\ind'}\mathbb E\big[\Corr_{\ind,\ind'}(y)\big] =  \Var(t_1).\]

The computation of the variance of $t_1$ is obtained by a computation similar  to that in second proof of Proposition \ref{WasDistLin} (page \pageref{Proof2}); we get
\begin{align*}
\Var(t_1) = & \sum_{m=k_0+1}^k \Var_{x_m\in[1/2-\lambda_{\fat^{m-1}(\ind)},1/2]}\left(x_m \widetilde\lambda_{\fat^{m}(\ind)}\right)\\
  = & \frac{1}{12} \sum_{m=k_0+1}^{k} \widetilde\lambda_{\fat^{m}(\ind)}^2 \lambda_{\fat^{m-1}(\ind)}^2 = \frac{1}{12} \sum_{m=k_0}^{k-1} \widetilde\lambda_{\fat^{m}(\ind)}^2.
\end{align*}
Finally,
\[\mathbb E\left[\Corr_{\ind,\ind'}(y)\right] =
\frac{1}{12} \sum_{m=k_0}^{k-1} \frac{\widetilde\lambda_{\fat^{m}(\ind)}^2}{\widetilde\lambda_{\ind} \widetilde\lambda_{\ind'}}.\]
\end{proof}

\begin{rem}
If in the formula we replace all the $\lambda_\ind$'s by $d$, we get
\[\Disc^2 = \frac{d^{k+1}-1}{12(d-1)}.\]
Without the correlations, it gives
\[\widetilde\Disc^2 = \frac{d^{k+1}-d^{-k-1}}{12(d-d^{-1})}.\]
the ratio between both tends to $1+1/d$ when $k$ goes to $+\infty$.
\end{rem}

We end this section by a quantitative version of Proposition \ref{WasDistTree}. For $E\subset \Z$, we will denote
\[D_R^+(E) = \sup_{x\in\R} \frac{\card\big(E\cap [x-R,x+R]\big)}{\card\big(\Z\cap [x-R,x+R]\big)}.\]

\begin{add}\label{AddBeurk}
For every $\ell', c\in\N$, there exists a locally finite union of positive codimension submanifolds $V_q$ of $]1,+\infty[^{\card I_k}$ such that for every $\eta'>0$, there exists a radius $R_0>0$ such that if $(\widetilde\lambda_{\ind})_{\ind\in I_k }$ satisfies $d\big((\widetilde\lambda_{\ind})_{\ind\in I_k }, V_q\big)>\eta'$ for every $q$, then for every $R\ge R_0$, and every family $\mathbf v = (v_\ind)_{\ind\in I_k}$ of real numbers, we have\footnote{The map $\widehat{\ell+v}$ is the discretization of the affine map $\ell+ v$.}
\begin{enumerate}[(i)]
\item\label{item1AddBeurk} there is a subset of $\Z$ with bounded gaps, made of points which are images of exactly $d^k$ points, each of them having all its roundoff errors close to 0: for every $y\in\Z$, there exists $y'\in\Z$ with $|y-y'|\le R_0$ and for every $\ind\in \llbracket 1,d\rrbracket^k$, a point $x_\ind$ such that 
\[\widehat{\Ell+\mathbf v}^k\big((x_\ind)_{\ind\in \llbracket 1,d\rrbracket^k}\big) = \{y'\}\]
and for every $\ind\in \llbracket 1,d\rrbracket^k$ and every $m\le k$, we have
\[|e_{x_\ind,\ind}^m| \le \frac{1}{\ell'},\]
where $e_{x,\ind}^m$ is the the error made at the $m$-th iteration of the point $x$, applying the discretizations of the maps $\ell_{\fat^j(\ind)} + v_{\fat^j(\ind)}$.
\item\label{item2AddBeurk} for any $y'$ like in item \eqref{item1AddBeurk}, the mean of the cumulated difference starting from $y'$ is almost zero (for a set $A\subset \R$, we denote $A-y'$ the translation of $A$ of $-y'$):
\[\left|\frac{1}{R}\int_0^R c\delta_x\left(\widehat{\Ell+\mathbf v}^k\big(\Z^{d^k}\big) - y'\, , \, \widetilde\lambda_{tot}^{-1} \Leb\right) \ud x\right| \le\frac{1}{\ell'}.\]
\item\label{item3AddBeurk} for any $y'$ like in item \eqref{item1AddBeurk}, the Cramér distance starting from $y'$ is almost the same as the one starting from 0:
\[
\left|\Disc^2_R \left(\widehat{\Ell+\mathbf v}^k\big(\Z^{d^k}\big) - y'\, , \, \widetilde\lambda_{tot}^{-1} \Leb\right)
 - \frac{1}{12\widetilde\lambda_{tot}} - \frac{1}{12} \sum_{\ind,\ind'\in \llbracket 1,d\rrbracket^k} \sum_{m=k_0(\ind,\ind')}^{k-1} \frac{\widetilde\lambda_{\fat^m(\ind)}\widetilde\lambda_{\fat^m(\ind')}}{\widetilde\lambda_{\ind}\widetilde\lambda_{\ind'}}\right| < \frac{1}{\ell'};
\]
\item\label{item4AddBeurk} there is only a small proportion of the points of the image sets which are obtained by discretizing points close to $\Z+1/2$: for every $m\le k$ and every $\ind\in \llbracket 1,d\rrbracket^k$, we have
\[D_R^+\left\{x\in \big(\ell_{\fat^m(\ind)} + v_{\fat^m(\ind)}\big)(\Z)\ \middle\vert\ d\big(x,\Z+\frac12 \big) < \frac{1}{3c\ell'} \right\} < \frac{1}{{c}\ell'};\]
\end{enumerate}
\end{add}

\begin{proof}[Sketch of proof of Addendum \ref{AddBeurk}]
The moral of this addendum is that if a collection of numbers $x_1,\dots,x_k$ is ``almost $\Q$-independent'', then the rotation of vector $x_1,\dots,x_k$ in $\T^k$ is ``ergodic up to $\varep$''.

In our particular case, if the collection $(\widetilde\lambda_{\ind})_{\ind\in I_k }$ does not satisfy any linear dependence relation with small integer coefficients, then 
for any $\ind,\ind'$, the image of the action of $\Z$ by $y\mapsto (0^{k+k_0-1},y)$ on $\R^{k+k_0}/ \widetilde M_{\ind,\ind'}\Z^{k+k_0}$ is ``uniformly distributed up to $\varep$''. This implies that the events $t_1$, $t_2$ and $t_2'$ are ``almost independent'' and that the variance of $t_1$ is almost equal to the formula of the proof of Proposition \ref{WasDistTree}.

These arguments are formalized by the following improvement of Weyl's criterion:

\begin{lemme}[Weyl]\label{Weyl}
Let $\dist$ be a distance generating the weak-* topology on $\Prb$ the space of Borel probability measures on $\T^n$. Then, for every $\varep>0$, there exists a locally finite family of affine hyperplanes $H_i \subset \R^n$, such that for every $\eta>0$, there exists $M_0\in\N$, such that for every $\boldsymbol\lambda \in \R^n$ satisfying $d(\boldsymbol\lambda,H_q)>\eta$ for every $q$, and for every $M\ge M_0$, we have
\[\dist \left(\frac{1}{M}\sum_{m=0}^{M-1} \bar\delta_{m\boldsymbol\lambda }\,,\ \Leb_{\R^n/\Z^n}\right) < \varep,\]
where $\bar\delta_x$ is the Dirac measure of the projection of $x$ on $\R^n/\Z^n$.
\end{lemme}

For a proof of this lemma, see \cite{Gui15a}.
\end{proof}


\section{A formula for the Cramér distance of $C^r$-generic expanding maps}\label{SecExpand}

In this section we prove Theorem~\ref{MainTheo}, by starting with the more explicit statement Theorem~\ref{ThWasDist}.

\subsection{First formula}

As a first step towards the proof of Theorem \ref{MainTheo}, one gets a first formula for the Cramér distance. Recall that for any $1\le r \le +\infty$, we denote $\mathcal D^r(\Sp^1)$ the set of $C^r$ expanding maps, and that we denote $d\ge 2$ the degree of $f\in\D^r(\Sp^1)$.

\begin{theoreme}\label{ThWasDist}
For any $1\le r \le +\infty$, if $f$ is a generic element of $\mathcal D^r(\Sp^1)$, then for any $k\in\N$, one has
\begin{multline}
\lim_{N\to +\infty} N  \Disc \big(f_N^k(E_N), \Ll^k(\Leb) \big)\\
\label{EqIntMieux}  = \left(\frac{1}{12} + \frac{1}{12} \int_{\Sp^1}\!\sum_{x,x'\in f^{-k}(y)} \sum_{m=k_0(x,x')}^{k-1}\!\! \frac{1}{(f^m)'(x) (f^m)'(x')} \ud y\right)^{1/2}.
\end{multline}
\end{theoreme}

Before coming to the proof, let us first make a few comments. First, note that when $f(x)=2x$, the right part of \eqref{EqIntMieux} becomes $2^{(k-3)/2}$. It gives an explicit approximation of the asymptotics for generic maps very close to $x\mapsto 2x$.

For $y \in \Sp^1$, let $H_y$ be the difference between the cumulative distribution functions of respectively:
\begin{itemize}
\item The uniform measure on $f_N^k(E_N)$, and
\item $\Ll^k(\Leb)$,
\end{itemize}
seen as measures on the fundamental domain $[y,y+1]$ of $\Sp^1$. By the definition of the Cramér distance (Equation \eqref{DefDisc}), 
\begin{equation}\label{EqIntMieux2}
\Disc \big(f_N^k(E_N), \Ll^k(\Leb) \big)^2 = \int_0^1 \left(H_y(x)-\Big(\int_0^1 H_y\Big) \right)^2 \ud x
\end{equation}

\begin{lemme}\label{LemHpFacil}
There exists a constant $B = B(f,k)>0$ such that for any $y$, one has $\|H_y\|_\infty \le B/N$.
\end{lemme}

\begin{proof}[Proof of Lemma \ref{LemHpFacil}]
The global roundoff error of each point $x\in E_N$ satisfies
\[\left| f_N^k(x) - f^k(x) \right|\le \frac{A}{N} \qquad \text{where} \qquad A = \frac{\|f'\|_\infty^k}{2(\|f'\|_\infty-1)}.\]
Thus, for any $y$, one has $\|H_y\|_\infty \le 2\frac{A}{N} \cdot \frac{d}{\min_{\Sp^1}f'}$.
\end{proof}

Hence, the good scale for the Cramér distance $\Disc \big(f_N^k(E_N), \Ll^k(\Leb) \big)$ is at most $1/N$. Theorem \ref{ThWasDist} ensures that this is exactly $1/N$, as it shows that $N\Disc \big(f_N^k(E_N), \Ll^k(\Leb) \big)$ converges towards a positive number when $N$ goes to infinity.
\bigskip

The proof of Theorem~\ref{ThWasDist} is mainly based on Proposition \ref{WasDistTree} (more precisely, Addendum \ref{AddBeurk}), which treats the linear corresponding case. By applying arguments of \cite{Gui15a}, one gets the following property:

\begin{prop}\label{propQuasiFormul}
Let $r\ge 1$ and $f$ a generic element of $\mathcal D^r(\T^n)$. Then for any $k\in\N$, any $\varep>0$, and any $N\in\N$ large enough, there exists a finite collection $I_p\subset \Sp^1$ of pairwise disjoint segments such that:
\begin{enumerate}[1)]
 \item\label{pt1} each segment $I_p$ has length smaller than $\varep$, and the union of the segments $I_p$ has Lebesgue measure bigger than $1-\varep$;
 \item\label{pt2} the left endpoint $y_p$ of each segment $I_p$ is an element of $E_N$;
 \item\label{pt3} each point $y_p$ is the image of $d^k$ points $x_{\ind,p}\in E_N$ (the maximal possible number) by $f_N^k$, and the roundoff error vector in time $k$ of each point $x_{\ind,p}$ is $\varep$-small;
 \item\label{pt4} for each $p$, the Cramér distance distribution restricted to the segment $I_p$ and starting from the point $y_p$ is $\varep$-close to the Cramér distance distribution associated to the preimage tree of $y_p$.
\end{enumerate}
\end{prop}

In particular, point (\ref{pt4}) ensures that on each segment $[y_p,y_p+R/N]\subset I_p$ such that $R\ge R_0$, the mean of the map $c\delta$ is $\varep$-close to 0, and its variance (which corresponds to the $L^2$ Cramér distance $\Disc$ associated to the map $f$) is $\varep$-close to the Cramér distance $\Disc$ associated to the preimage tree of $y_p$.

\begin{proof}[Proof of Proposition \ref{propQuasiFormul}]
We simply apply the arguments of the proof of Theorem 33 of \cite{Gui15a}, by replacing Lemma 34 of \cite{Gui15a} by Addendum \ref{AddBeurk}. In particular, this proof tells us that Thom's transversality theorem implies the existence of a family $[y_p',z_p]$ of segments of length $\gg R_0 N$ (where $R_0$ is given by Addendum \ref{AddBeurk}), with $y_p',z_p\in 1/N \Z$. For any $p$, we apply Addendum \ref{AddBeurk} to the point $y=y_p'$, the derivatives given by the preimage tree starting from $y_p'$, and the preimage set of $y_p$ as the vector $\mathbf v$.
This gives us a point $y'=y_p$ (by item \eqref{item1AddBeurk}), and allows to define the segments $I_p = [y_p,z_p]$.

Point \eqref{pt1} comes from the proof of Theorem 33 of \cite{Gui15a} and Points (\ref{pt2}) and (\ref{pt3}) come from point \eqref{item1AddBeurk} of Addendum \ref{AddBeurk}. For itself, point \eqref{pt4} comes from an application of Taylor formula, the linear formulation of items item \eqref{item2AddBeurk} and item \eqref{item3AddBeurk} of Addendum \ref{AddBeurk}, and the error estimate for nonlinearities of item \eqref{item4AddBeurk} of Addendum \ref{AddBeurk} (see the proof of Theorem 33 of \cite{Gui15a} for more details).
\end{proof}

\begin{proof}[Proof of Theorem \ref{ThWasDist}]
In this proof we denote $H_p = H_{y_p}$, where the $y_p$ are given by Proposition \ref{propQuasiFormul}.

Points (\ref{pt2}) and (\ref{pt3}) of Proposition \ref{propQuasiFormul} ensure that for $p\neq p'$, the cumulative distribution functions $H_p$ and $H_{p'}$ (seen as functions of $\Sp^1$) are close: reasoning as in Lemma \ref{LemHpFacil}, one gets
\begin{equation}\label{EqHpClose}
\|H_p-H_{p'}\|_\infty \le \varep \frac{B}{N}.
\end{equation}

By Taylor formula, if $\varep$ is small enough, in restriction to the interval $I_p$, the measure $\Ll^k(\Leb)$ is close to $\widetilde\lambda_{tot}^{-1} \Leb$, where $\widetilde\lambda_{tot}$ denotes the multiplier associated to the preimage tree at $y_p$. Combined with point (\ref{pt4}), this fact implies that the mean of $H_p$ restricted to $I_p$ is small:
\[\left|\frac{1}{|I_p|}\int_{I_p} H_p\right| \le \frac{\varep}{N}.\]
Hence, fixing $p=0$, and using the fact that all the $H_p$'s are close, one gets that the mean of $N H_0$ is small in restriction to the union of the $I_p$.

Then, using point (\ref{pt1}) together with Lemma \ref{LemHpFacil}, we deduce that the mean of $N H_0$ is small: there is a constant $C=C(f,k)>0$ such that 
\[\left|\int_{\Sp^1} H_0 \right| \le \varep \frac{C}{N}.\]
This fact, together with Equation \eqref{EqIntMieux2}, implies that
\[\left|\Disc \big(f_N^k(E_N), \Ll^k(\Leb) \big)^2 - \int_{\Sp^1} H_0^2 \right| \le \frac{\varep}{N^2}.\]
Using Equation \eqref{EqHpClose} and Lemma \ref{LemHpFacil} again, we deduce that
\[\left|\Disc \big(f_N^k(E_N), \Ll^k(\Leb) \big)^2 - \sum_p \int_{I_p} H_p^2 \right| \le \frac{\varep}{N^2}.\]

We now use point (\ref{pt4}), which ensures that for any $p$,
\[ \left|N^2\int_{I_p} H_p^2 - \left(\frac{1}{12} + \frac{1}{12} \int_{I_p}\!\sum_{x,x'\in f^{-k}(y)} \sum_{m=k_0(x,x')}^{k-1}\!\! \frac{1}{(f^m)'(x) (f^m)'(x')} \ud y\right) \right| \le \varep |I_p|.\]
Combined with the previous estimation, this gives the theorem.
\end{proof}

\subsection{Proof of theorem \ref{MainTheo}}


\begin{proof}[Proof of Theorem \ref{MainTheo}]
Recall that in Theorem \ref{ThWasDist} we have, for any $x_0$ that is an $f^r$-preimage of $y$ for some $r$:
\[
   \widetilde{\lambda}_{x_0} = f'(x_0) \cdot f'(f(x_0)) \cdot \dots \cdot f'(f^{r-1}(x_0)) =
   \frac{\ud}{\ud x}(f^r(x))|_{x=x_0} = Df^{r}(x_0),
 \]
 and for each $x,x'\in f^{-k}(y)$ we have $k_0(x,x')$ defined to be the smallest $m$ such
 that $f^m(x)=f^m(x')$.

Therefore the integral of Theorem \ref{ThWasDist} can also be written as
\begin{equation}\label{EqIntMelhorouAindaMais}
\int_0^1\sum_{x,x'\in f^{-k}(y)} \sum_{m=k_0(x,x')}^{k-1} \frac{1}{Df^{m}(x) \cdot Df^{m}(x')} \ud y.
\end{equation}
Putting $n+e=k$ (for each pair $x,x'$, being the $f$-images \emph{n}ot equal for $n$ steps, then \emph{e}qual for $e$ steps), we can split this integral in sums of the form
\[
  \sum_{\substack{n+e=k\\n,e\geq 0}}  \int_0^1 \sum_{z\in f^{-e}(y)}
  \sum_{\substack{x,x'\in f^{-n}(z)\\f^{n-1}(x)\neq f^{n-1}(x')}} \sum_{m=n}^{k-1}
  \frac{1}{Df^{m}(x) Df^{m}(x')} \ud y.
\]
When we fix a certain $0\leq m \leq k-1$, the integral is
\begin{multline*}
   \int_0^1 \sum_{\substack{n\leq m\\e=k-n}} \sum_{z\in f^{-e}(y)}
  \sum_{\substack{x,x'\in f^{-n}(z)\\f^{n-1}(x)\neq f^{n-1}(x')}} 
  \frac{1}{Df^{m}(x) Df^{m}(x')} \ud y \\
  = \int_0^1 \sum_{w\in f^{-(k-m)}(y)}\ 
  \sum_{\substack{x,x'\in f^{-m}(w)}} 
  \frac{1}{Df^{m}(x) Df^{m}(x')} \ud y
\end{multline*}
as it could be also deduced directly from Equation \eqref{EqIntMelhorouAindaMais} fixing
$m$.

Since $y = f^{k-m}(w)$ and therefore $\ud y=Df^{k-m}(w)\ud w$, dividing the domain into the intervals where $f^{k-m}$ is injective, and then putting again everything together, we can change variable and obtain
\begin{equation} \label{IntSumW1}
  \int_0^1 \sum_{\substack{x,x'\in f^{-m}(w)}} 
  \frac{Df^{k-m}(w)}{Df^{m}(x) Df^{m}(x')} \ud w.
\end{equation}

Let us consider now the map $f\times f:[0,1]^2\rightarrow [0,1]^2$, and let
$L_{f\times f}$ be its transfer operator. Given an observable
$H(w,w'):[0,1] \rightarrow \R$, the $m$-th power of $L_{f\times f}$ can be
computed on $H$ as
\[
   \big(L_{f\times f}^m H\big)(w,w') = \sum_{\substack{f^{m}(x)=w\\f^{m}(x')=w'}} \frac{H(x,x')}{Df^{m}(x) Df^{m}(x')},
\]
 considering that the Jacobian determinant of $(f\times f)^m$ at $(x,x')$ is $Df^{m}(x) Df^{m}(x')$. Notice also that if $H(x,x')=h_1(x)\cdot h_2(x')$ we have that
\[
    \big(L^m_{f\times f}H\big)(w,w') = (L^m_fh_1)(w) \cdot (L^m_fh_2)(w').
\]

 Therefore, the integral of \eqref{IntSumW1} is also the integral on the diagonal
 \[
   \Delta = \{(w,w) :
 w\in[0,1]\}
 \]
 of $L_{f\times f}^m H$ for any observable $H$ such that $H(x,x')=Df^{k-m}(w)$ whenever $f^m(x)=w=f^m(x')$. In our case can take for instance
 \[
     H(x,x') = \sqrt{Df^{k-m}(f^{m}x) \cdot Df^{k-m}(f^{m}x')},
\]
and in the end the integral amounts to the integral $L^m_{f\times f}(H)$ along $\Delta$.

Taking $h(w)=\sqrt{Df^{k-m}(f^{m}w)}$ we have $H(w,w')=h(w)\cdot h(w')$, therefore
\begin{align*}
 \int_0^1 \sum_{\substack{x,x'\in f^{-m}(w)}} 
    \frac{Df^{k-m}(w)}{Df^{m}(x) Df^{m}(x')} \ud w
  = &\int_0^1 (L^m_{f\times f}H)(w,w) \ud w \\
  = &\int_0^1 (L^m_{f}h)(w) \, (L^m_{f}h)(w) \ud w \\
  = &\int_0^1 (L^m_{f}h)(w)^2 \ud w.
\end{align*}

Let $g(x) = \sqrt{Df^{k-m}(x)}$, so that $h(x) = g(f^{m}(x))$ for short. Let us recall the formula for the transfer operator applied to $g(f^{m}(x))$:
\begin{align*}
  (L_f^m h)(w) = &\sum_{x\in f^{-m}(w)} \frac{g(f^{m}(x))}{Df^{m}(x)} \\
  = & \sum_{x\in f^{-m}(w)} \frac{g(w)}{Df^{m}(x)} \\
  = & g(w) \sum_{x\in f^{-m}(w)} \frac{1}{Df^{m}(x)} \\
  = & g(w) \, (L_f1)(w).
\end{align*}
Therefore, our integral can be written as
\begin{align*}
\int_0^1 (L^m_{f}h)(w)^2 \ud w 
  = & \int_0^1 \left[ g(w) \, (L_f^m 1)(w) \right]^2 \ud w \\
  = & \int_0^1 Df^{k-m}(w) \, (L_f^m 1)(w)^2 \ud w \\
  = & \big\langle Df^{k-m}, (L_f^m 1)^2 \big\rangle.
\end{align*}

Taking the sum over $m=0,1,2,\dots, k-1$, we have proved Theorem \ref{MainTheo}.
\end{proof}


\begin{rem}
Stating from Equation \eqref{IntSumW1}, we can evaluate the case where the sum is restricted to $x=x'$. Changing variable $w=f^{m}(x)$ have
\begin{align*}
\int_0^1 \sum_{\substack{x\in f^{-m}(w)}} \frac{Df^{k-m}(w)}{\left[Df^{m}(x)\right]^2} \ud w
     = & \int_0^1 \frac{Df^{k-m}(f^{m}(x))}{\left[Df^{m}(x)\right]^2} Df^{m}(x)\ud x \\
     = & \int_0^1 \frac{Df^{k-m}(f^{m}(x))}{Df^{m}(x)} \ud x \\
     = & \int_0^1 \frac{Df^{k}(x)}{\left[Df^{m}(x)\right]^2} \ud x.
\end{align*}
\end{rem}

\bibliographystyle{amsalpha}
\bibliography{../Biblio.bib}
\vfill

\end{document}